\RequirePackage{fix-cm}
\documentclass[smallextended]{svjour3}       % onecolumn (second format)
\usepackage{graphicx}
\usepackage{algorithm}
\usepackage{algpseudocode}
\usepackage{multirow}
\usepackage{tabulary}

\usepackage{amsfonts}
\usepackage{amsmath}
\usepackage{amssymb}

\usepackage{xcolor}
\usepackage{color}
\usepackage{soul}
\begin{document}

\title{Direct simple computation of middle surface between 3D point clouds and/or discrete surfaces by tracking sources in distance function calculation algorithms}
\thanks{This work was supported by grants APVV-19-0460, and VEGA 1/0436/20.
%\thanks{
%Grants or other notes
}

\titlerunning{Direct computation of middle surfaces}        % if too long for running head

\author{Balázs Kósa         \and
        Karol Mikula %etc.
}


\institute{B. Kósa \at
              %first address \\
              \email{kosa@math.sk}           %  \\
           \and
           K. Mikula \at
           \email{mikula@math.sk}\\
           %second address
           Department of Mathematics and Descriptive Geometry, Faculty of Civil Engineering, Slovak University of Technology in Bratislava, Slovakia
}


\maketitle

\begin{abstract}
In this paper, we introduce novel methods for computing middle surfaces between various 3D data sets such as point clouds and/or discrete surfaces. Traditionally the middle surface is obtained by detecting singularities in computed distance function such as ridges, triple junctions, etc. It requires to compute second order differential characteristics and also some kinds of heuristics must be applied. Opposite to that, we determine the middle surface just from computing the distance function itself which is a fast and simple approach. We present and compare the results of the fast sweeping method, the vector distance transform algorithm, the fast marching method, and the Dijkstra-Pythagoras method in finding the middle surface between 3D data sets.
\end{abstract}

\section{Introduction} \label{intro}

\par Finding an optimal middle surface for a data set is a crucial task in many applications such as computational geometry, surface representation and reconstruction, image processing and computer vision or mesh generation.
In optimal mesh generation \cite{Persson_phdthesis} for example, the information about the middle surface can be used to densify or coarsen the computational grid in the computational domain. For this reason, having an efficient method that fulfills such needs is very important. Very often the middle surface for which algorithms are seeking is a middle axis of a closed curve or a surface, see e.g. \cite{KIMMEL1995382,Siddiqi_Hamilton,Rumpf_ContSkelet}. Such algorithms can be complicated because they utilize second order derivatives of the computed distance function in order to detect its ridges, junctions and other singularities which often requires some kinds of heuristics, see also \cite{Persson_phdthesis}. In cases where we can distinguish individual separate or labelled shapes between which we want to find the middle surface, a much more straightforward approach can be derived. We show how algorithms designed for distance function calculation can be adjusted and utilized in these cases to obtain the middle surface already during the computations of the distance function itself. Opposite to methods that utilize second order derivatives of the computed distance function, we only adjust the distance function calculation algorithms. This makes our methods simple, efficient and easy to implement. 
\par A distance function to an object is a useful tool in a variety of disciplines. For this reason, over the years many algorithms have been developed which were optimized to obtain the most accurate result as fast as possible, see e.g. \cite{Jones_Baerentzen_Sramek}. We provide a short description of four such algorithms and show how they can be implemented to calculate the distance function on a uniform voxel grid for 3D objects represented either by a point clouds or triangulated surfaces. To compare the algorithms, we applied them to several data sets and measured their accuracy and speed.
\par After providing a sufficient explanation of the methods with detailed pseudo-codes for each of them, we describe how we use them to find the middle surface. We will see that all it needs is a few natural changes in the implementation to achieve this goal.  We test our approaches on several experiments and present the results subsequently.

%%%%%%%%%%%%%%%%%%%%%%%%%%%%%%%%%%%%%%%%%%%%%%%%%%%%%%%%%%%%%%%%%%%%%%Numerical methods%%%%%%%%%%%%%%%%%%%%%%%%%%%%%%%%%%%%%%%%%%%%%%%%%%%%%%%%%%%%%%%%%%%%%%

\section{Numerical methods} \label{Numerical_methods}

In computational mathematics, the notion of distance function is used for the result of distance computation. In this section, we will discuss common numerical methods used for this task. Following \cite{Smisek_dissertation} the presented methods are classified according to following two criteria:
\begin{enumerate}
	\item Distance definition: The distance function can be calculated as a solution of the so-called eikonal equation or by the Euclidean distance computation.
	\item Voxel visit order strategy: We will analyze sweeping and wavefront methods.
\end{enumerate}
Our goal is to demonstrate how methods falling under these categories can be used to find the middle surface between two or more input data sets. We will discuss and analyze these methods: the fast sweeping method (FSM) \cite{Zhao}, the vector distance transform (VDT) algorithm \cite{Danielsson}, the fast marching method (FMM) \cite{Sethian_FMM} and the Dijkstra-Pythagoras (DP) method \cite{Smisek_dissertation}. Table \ref{tab:distance_function_methods} shows the classification of the four studied methods. 

\begin{table}[h]
	\centering
	\begin{tabular}{|c|c|c|}
		\hline
		\multicolumn{1}{|l|}{\textbf{}} & \textbf{Sweeping}                  & \textbf{Wavefront}            \\ \hline
		\textbf{Eikonal equation}       & \textbf{Fast sweeping method}      & \textbf{Fast marching method} \\ \hline
		\textbf{Euclidean distance}     & \textbf{Vector distance transform} & \textbf{Dijkstra-Pythagoras}  \\ \hline
	\end{tabular}
	\caption{Distance function computation methods category overview. Rows represent the distance definition. Columns represent the voxel visit order strategy.}
	\label{tab:distance_function_methods}
\end{table}

%%%%%%%%%%%%%%%%%%%%%%%%%%%%%%%%%%%%%%%%%%%%%%%%%%%%%%%%%Explanation of classification%%%%%%%%%%%%%%%%%%%%%%%%%%%%%%%%%%%%%%%%%%%%%%%%%%%%%%%%%%%%%%%%%%%%

\subsection{Basic definitions} \label{Basic_definitions}

The distance function will be calculated on the computational domain $\Omega$, $\Omega \subseteq \mathbb{R}^{n}$. A data set $\Omega_{0}$, to which we want to compute distance function $d$, will be a subset of $\Omega$, $\Omega_{0} \subseteq \Omega$. In this paper, we work with 3D objects so we limit the dimension of $\Omega$ to $n=3$. With this notation of the domain we can define the distance function as $d:\Omega \rightarrow \mathbb{R}$. On subset $\Omega_{0}$ the distance should be $0$, thus we get the boundary condition
\begin{equation} \label{eq:distance_ini}
	d\left(x\right)=0, \; x \in \Omega_{0}\subseteq \Omega. 
\end{equation}
Then the task is to calculate $d\left(x\right), \; x \in \Omega\setminus \Omega_{0}$.

\subsubsection{Distance definition} \label{Distance_definition}

For the numerical methods, the computational domain $\Omega$ will be discretized into a finite number of voxels with edge size $h$. The number of voxels will be denoted as $N_{i}$ along the $x$ axis, $N_{j}$ along the $y$ axis and $N_{k}$ along the $z$ axis. In the obtained computational grid, the function $d$ will be calculated at the center of every voxel, the so-called grid points.
\par The \textbf{eikonal equation} is given by
\begin{equation} \label{eq:eikonal}
	\left |\nabla d \left ( x \right )  \right | = 1 \quad x \in \Omega.
\end{equation}
This equation will be coupled with the boundary condition \eqref{eq:distance_ini}. For the discretization of \eqref{eq:eikonal}, we denote grid points of $\Omega$ by $x_{i,j,k}$ and the numerical solution of the distance function at $x_{i,j,k}$ as $d_{i,j,k}$. The discretization of \eqref{eq:eikonal} at interior grid points is done according to the Godunov upwind difference scheme \cite{Rouy_Tourin}:
\begin{equation} \label{godunov}
	\begin{split}
		\left [ \left ( d_{i,j,k}-d_{x\, min} \right )^{+} \right ]^{2} &+\left [ \left ( d_{i,j,k}-d_{y\, min} \right )^{+} \right ]^{2} + \left [ \left ( d_{i,j,k}-d_{z\, min} \right )^{+} \right ]^{2} = h^{2},\\
		i=1,...&,I-1, \; j=1,...,J-1, \; k=1,...,K-1,\\
		&d_{x\, min} =min\left ( d_{i,j-1,k},\: d_{i,j+1,k} \right ),\\
		&d_{y\, min} =min\left ( d_{i-1,j,k},\: d_{i+1,j,k} \right ),\\
		&d_{z\, min} =min\left ( d_{i,j,k-1},\: d_{i,j,k+1} \right ),
	\end{split}
\end{equation}
$$\left (x  \right )^{+}=\left\{\begin{matrix}
	x,\; x>0\\ 
	0,\; x\leq0
\end{matrix}\right. \rm{.}$$
At the boundary of $\Omega$ we use one sided difference. This enforces that the solution at every voxel center is defined by the smaller values of neighboring grid points. Eikonal-based methods calculate the distance function by applying the described numerical scheme (\ref{godunov}).
\par \textbf{Euclidean distance} between two points will be defined according to the Pythagoras’ theorem. For $a = \left ( a_{x}, a_{y}, a_{z} \right ) \in \Omega, b = \left ( b_{x}, b_{y}, b_{z} \right ) \in \Omega$ we define
\begin{equation} \label{eq:pythagoras}
	d\left ( a, b \right ) = \sqrt{\left ( a_{x} -b_{x} \right )^{2} + \left ( a_{y} -b_{y} \right )^{2} + \left ( a_{z} -b_{z} \right )^{2}}. 
\end{equation}

\subsubsection{Voxel visit order strategy} \label{Voxel_visit_order_strategy}

For algorithms with the \textbf{sweeping} approach, Gauss-Seidel iterations with alternating sweeping orderings are used. This allows the methods to pass through the voxels multiple times. For three dimensions we sweep the computational domain with eight alternating orderings:
\begin{equation*} %\label{sweeps}
	\begin{split}
		& 1.\: i=1:N_{i},\: j=1:N_{j},\: k=1:N_{k}; \;\; 2.\: i=1:N_{i},\: j=1:N_{j},\: k=N_{k}:1; \\  
		& 3.\: i=1:N_{i},\: j=N_{j}:1,\: k=1:N_{k}; \;\; 4.\: i=1:N_{i},\: j=N_{j}:1,\: k=N_{k}:1; \\
		& 5.\: i=N_{i}:1,\: j=1:N_{j},\: k=1:N_{k}; \;\; 6.\: i=N_{i}:1,\: j=1:N_{j},\: k=N_{k}:1; \\
		& 7.\: i=N_{i}:1,\: j=N_{j}:1,\: k=1:N_{k}; \;\; 8.\: i=N_{i}:1,\: j=N_{j}:1,\: k=N_{k}:1.
	\end{split}
	%\begin{split}
	%	& 1.\: i=1:N_{i},\: j=1:N_{j},\: k=1:N_{k};\\ 
	%	& 2.\: i=1:N_{i},\: j=1:N_{j},\: k=N_{k}:1;\\  
	%	& 3.\: i=1:N_{i},\: j=N_{j}:1,\: k=1:N_{k};\\
	%	& 4.\: i=1:N_{i},\: j=N_{j}:1,\: k=N_{k}:1;\\
	%	& 5.\: i=N_{i}:1,\: j=1:N_{j},\: k=1:N_{k};\\
	%	& 6.\: i=N_{i}:1,\: j=1:N_{j},\: k=N_{k}:1;\\
	%	& 7.\: i=N_{i}:1,\: j=N_{j}:1,\: k=1:N_{k};\\   
	%	& 8.\: i=N_{i}:1,\: j=N_{j}:1,\: k=N_{k}:1.
	%\end{split}
\end{equation*}
To work with these sweeps in the following sections we will define the following sets:
\begin{equation} \label{sweeps}
	\begin{split}
	& i_{sweep} = \left \{ \left \{ 0,N_{i}-1,1 \right \}, \left \{ 0,N_{i}-1,1 \right \}, \left \{ 0,N_{i}-1,1 \right \}, \left \{ 0,N_{i}-1,1 \right \}, \right.\\
	& \left. \left \{ N_{i}-1,0,-1 \right \}, \left \{ N_{i}-1,0,-1 \right \}, \left \{ N_{i}-1,0,-1 \right \}, \left \{ N_{i}-1,0,-1 \right \} \right \} \\
	& j_{sweep} = \left \{ \left \{ 0, N_{j} - 1, 1 \right \}, \left \{ 0, N_{j} - 1,  1 \right \}, \left \{ N_{j} - 1, 0, -1 \right \}, \left \{ N_{j} - 1, 0, -1 \right \}, \right.\\ 
	& \left. \left \{ 0, N_{j} - 1,  1 \right \}, \left \{ 0, N_{j} - 1,  1 \right \}, \left \{ N_{j} - 1, 0, -1 \right \}, \left \{ N_{j} - 1, 0, -1 \right \} \right \} \\
	& k_{sweep} = \left \{ \left \{ 0, N_{k} - 1, 1 \right \}, \left \{ N_{k} - 1, 0, -1 \right \}, \left \{ 0, N_{k} - 1,  1 \right \}, \left \{ N_{k} - 1, 0, -1 \right \}, \right. \\
	& \left. \left \{ 0, N_{k} - 1, 1 \right \}, \left \{ N_{k} - 1, 0, -1 \right \}, \left \{ 0, N_{k} - 1, 1 \right \}, \left \{ N_{k} - 1, 0, -1 \right \} \right \}.
	\end{split}
\end{equation}
The different algorithms analyze a certain set of neighboring voxels in every iteration. This can be the set of 6 closest neighbors
\begin{equation} \label{set:6_neighbors}
P^{1}=\left\{\left(r,s,t\right);r,s,t\in\left\{-1,0,1\right\};\left|r\right|+\left|s\right|+\left|t\right|=1\right\}
\end{equation}
or the set including also the diagonal voxels, the set of all 26 neighbors 
\begin{equation} \label{set:26_neighbors}
	P^{2}=\left\{\left(r,s,t\right);r,s,t\in\left\{-1,0,1\right\};\left|r\right|+\left|s\right|+\left|t\right|=c;c \in \left\{ 1,2,3 \right\}\right\}.
\end{equation}
\par In the \textbf{wavefront methods} at every grid point, we assign the final value already in the first pass. To ensure this, the algorithms have to be set up in a way that every voxel is visited in the correct order, starting with the voxel nearest to $\Omega_{0}$ and ending with the furthest. For this, a data structure called \textit{min-priority-heap} \cite{Intro_algo} is utilized. In this structure whenever a change occurs the elements are rearranged so the element with the smallest value is on top. For wavefront algorithms at the beginning, we store all grid points that enforce the boundary condition in such a heap, with their distance value $d$ and their location in the grid. In every iteration, we can immediately obtain the grid point with the smallest value of $d\left(x\right)$. As the front moves on, new elements are added to the heap. For easy updates of distance values at gird points already saved in the heap, additional information about their location in the heap should be maintained.
\par In the next subsections, we will go through the implementation of the mentioned methods, so we will be able to describe how to change them for the task of computing the middle surface. To that goal, we start with the description of how to implement the initialization of the distance function to ensure the boundary condition \eqref{eq:distance_ini}.

%%%%%%%%%%%%%%%%%%%%%%%%%%%%%%%%%%%%%%%%%%%%%%%%%%%%%%%%%%%%%%%%%%%%%%Initialization%%%%%%%%%%%%%%%%%%%%%%%%%%%%%%%%%%%%%%%%%%%%%%%%%%%%%%%%%%%%%%%%%%%%%%

\subsection{Initialization} \label{Initialization_phase}
For every point $x$ of the input data set $\Omega_{0}$ the function $d$ should fulfill \eqref{eq:distance_ini}. When we implement a method for the calculation of $d$, we need to find a way to fulfill this condition. If point $x$ would coincide with the voxel center, in an array representing $d$ we could just set the values to $0$ for every such point $x$. Unfortunately, most of the time this is not the case. 
\par While working with point cloud data, to fulfill \eqref{eq:distance_ini}, we initialize the function $d$ as follows. We find the 8 nearest grid points to every point in the cloud and calculate the exact distance for these points from the corresponding point cloud element. The smallest possible distance will be saved at grid points when exploring subsequently all point cloud elements. These initialized values will be fixed in further calculations. Some of the algorithms described in the following sections use the cloud points as "sources" to calculate the distance function at other grid points. For this reason, in the initialization, we will keep track of this information as well. We can easily do this by setting the index of the source cloud point to the fixed grid points which will refer to the coordinates of the source. At other than fixed grid points, we set $d$ to a high enough number, which is bigger than the biggest possible distance in the grid. To simplify this, we can use $+\infty$, which for example when we implement the algorithm in C or C++ can be substituted by the maximum $double$ value.
\par In Alg. \ref{Initialization_pc}, we show how the described initialization can be easily implemented.

\clearpage

\begin{algorithm}
	\caption{Initialization of distance function to the point cloud data}
	\label{Initialization_pc}
	\begin{algorithmic}[1]
		\renewcommand{\algorithmicrequire}{\textbf{Input:}}
		\Require \texttt{Point cloud data:} \linebreak
		\texttt{$pc_{l}$ - $\left(x,y,z\right)$ coordinates of the $l$th point,} \linebreak
		\texttt{$N$ - number of points.}
		%\mbox{\texttt{~~~3D grid with voxel edge size $h$ and dimensions:}} \linebreak 
		\Require \texttt{3D grid with voxel edge size $h$ and dimensions $N_{i}, N_{j}, N_{k}$.}
		\renewcommand{\algorithmicrequire}{\textbf{Declaration:}}
		\Require \texttt{Arrays:} \linebreak 
		\texttt{$d_{i,j,k}$ - value of distance function at grid point $\left(i,j,k\right)$, } \linebreak
		\texttt{$c_{i,j,k}$ - $\left(x,y,z\right)$ coordinates of grid point $\left(i,j,k\right)$, } \linebreak
		\texttt{$f_{i,j,k}$ - determines if $d_{i,j,k}$ is fixed at $(i,j,k)$, } \linebreak
		\texttt{$s_{i,j,k}$ - source for $d_{i,j,k}$ calculation at $(i,j,k)$. } 
		\State \textbf{Set:} \texttt{$d_{i,j,k}$ to $+\infty$, $f_{i,j,k}$ to $false$, $s_{i,j,k}$ to $unknown$ }
		\State \textbf{Calculate:} $c_{i,j,k}$
		\For{$\left(l=0;l<N;l=l+1\right)$}
		\State $i_{first} = RoundDown\left(\left(pc_{l}.x - min\left(c_{i,j,k}.x\right)\right)/h\right)$
		\State $j_{first} = RoundDown\left(\left(pc_{l}.y - min\left(c_{i,j,k}.y\right)\right)/h\right)$
		\State $k_{first} = RoundDown\left(\left(pc_{l}.z - min\left(c_{i,j,k}.z\right)\right)/h\right)$
		\For{$\left(i=i_{first};i \leq i_{first}+1;i=i+1\right)$}
		\For{$\left(j=j_{first};j \leq j_{first}+1;j=j+1\right)$}
		\For{$\left(k=k_{first};k \leq k_{first}+1;k=k+1\right)$}
		\State $d_{new} = d\left(pc_{l},c_{i,j,k}\right)$	 \Comment{Calculated by \eqref{eq:pythagoras}.}
		\If{$d_{new}<d_{i,j,k}$}
		\State $d_{i,j,k}=d_{new}$ 
		\State $f_{i,j,k}=true$ 
		\State $s_{i,j,k}=pc_{l}$
		\EndIf
		\EndFor	
		\EndFor	
		\EndFor	
		\EndFor
	\end{algorithmic}
\end{algorithm}

%%%%%%%%%%%%%%%%%%%%%%%%%%%%%%%%%%%%%%%%%%%%%%%%%%%%%%%%%%%%%%%%%Fast sweeping method%%%%%%%%%%%%%%%%%%%%%%%%%%%%%%%%%%%%%%%%%%%%%%%%%%%%%%%%%%%%%%%%%%

\subsection{Fast sweeping method}
The fast sweeping method (FSM) \cite{Zhao} is an iterative algorithm with alternating sweeps \eqref{sweeps} used for the numerical solution of the Eikonal equation \eqref{godunov}. It can be applied in any number of dimensions for a rectangular computational grid. The value of $d\left ( x \right )$ at any grid point will never increase because an update rule is implemented by which the new value of the distance function is saved only if it is smaller than the current value. This enforces the correct value not to change at later iterations. 
\par \textcolor{black}{Let us denote in equation \eqref{godunov} the unknown as $x = d_{i,j,k}$ and the coefficients as $a_{1}=d_{x\, min}$, $a_{2}=d_{y\, min}$, $a_{3}=d_{z\, min}$. Then the unique solution, denoted by $\bar{x}$, to the equation
\begin{equation} \label{quadratic}
	\left [ \left ( x-a_{1} \right )^{+} \right ]^{2} +\left [ \left ( x-a_{2} \right )^{+}  \right ]^{2} + \left [ \left ( x-a_{3} \right )^{+}  \right ]^{2} = h^{2}
\end{equation}
can be found as follows. We order $a_{1},a_{2},a_{3}$ in increasing order. For generality we assume $a_{1}\leq a_{2}\leq a_{3}$. There is an integer $p,\,1 \leq p \leq 3$, such that $\bar{x}$ is the unique solution that satisfies
\begin{equation} \label{quadratic2}
	\left ( x-a_{1} \right )^{2} +\left ( x-a_{2} \right )^{2} + \left ( x-a_{3} \right )^{2} = h^{2} \quad and \quad a_{p}<\bar{x}<a_{p+1}
\end{equation}
To find $\bar{x}$ we start with $p=1$. If $\tilde{x}=a_{1}+h\leq a_{2}$ then $\bar{x}=\tilde{x}$. Otherwise we have to find the solution of the quadratic equation $$\left ( x-a_{1} \right )^{2} +\left ( x-a_{2} \right )^{2}= h^{2} $$ that satisfies $\tilde{x}>a_{2}$. We always take the maximum of the two solutions as our $\tilde{x}$. If $\tilde{x} \leq a_{3}$ then $\bar{x}=\tilde{x}$. If we still doesn't have a $\tilde{x}$ which satisfies all the conditions as the third step we compute the solution of the quadratic equation
$$\left ( x-a_{1} \right )^{2} +\left ( x-a_{2} \right )^{2} + \left ( x-a_{3} \right )^{2}= h^{2} $$ which will satisfy \eqref{quadratic2}.}
\par Only a finite amount of iterations is needed to obtain the solution, thus the complexity of the method is $O\left(N\right)$, where $N$ is the total number of grid points in the computational domain. This method is simple to implement, as it can be seen in the provided pseudo-code Alg. \ref{FSM_pseudo_code}. 

\begin{algorithm}
	\caption{Fast sweeping method}
	\label{FSM_pseudo_code}
	\begin{algorithmic}[1]
		\renewcommand{\algorithmicrequire}{\textbf{Input:}}
		\Require \texttt{From Alg.\ref{Initialization_pc}: 3D grid, $d_{i,j,k}$, $f_{i,j,k}$ }
		%\renewcommand{\algorithmicrequire}{\textbf{Definition:}}
		%\Require \texttt{Eight alternating orderings \eqref{sweeps}.}
		\For{$\left(l=0;l<8;l=l+1\right)$}
		\For{$\left(i=i_{sweep}\left[l,0\right];i \leq i_{sweep}\left[l,1\right];i=i+i_{sweep}\left[l,2\right]\right)$}
		\For{$\left(j=j_{sweep}\left[l,0\right];j \leq j_{sweep}\left[l,1\right];j=j+j_{sweep}\left[l,2\right]\right)$}
		\For{$\left(k=k_{sweep}\left[l,0\right];k \leq k_{sweep}\left[l,1\right];k=k+k_{sweep}\left[l,2\right]\right)$}
		\If{$f_{i,j,k}$ \texttt{is not true} }
		\State $a_{1} = min\left(d_{i+1,j,k}, d_{i-1,j,k}\right)$ 
		\State $a_{2} = min\left(d_{i,j+1,k}, d_{i,j-1,k}\right)$
		\State $a_{3} = min\left(d_{i,j,k+1}, d_{i,j,k-1}\right)$ \Comment{Use $+\infty$ if $\left(i,j,k\right)$ is out of bounds.}
		\State \texttt{Sort $\left \{ a_{1}, a_{2}, a_{3} \right \}$ from lowest to highest.}
		\State $d_{new}=a_{1}+h$
		\If{$d_{new}>a_{2}$}
		\State $d_{new} = \underset{x}{MaxSolution}\left( \left( x - a_{1}  \right)^{2} + \left( x - a_{2}  \right)^{2} = h^{2} \right)$
		\If{$d_{new}>a_{3}$}
		\State $d_{new} = \underset{x}{MaxSolution}\left( \left( x-a_{1}  \right)^{2} + \left( x-a_{2}  \right)^{2} + \left( x-a_{3}  \right)^{2} = h^{2} \right)$
		\EndIf
		\EndIf
		%\If{$d_{new}<d_{i,j,k}$}
		%\State $d_{i,j,k}=d_{new}$
		%\EndIf
		\State \textbf{if} \texttt{$d_{new}<d_{i,j,k}$} \textbf{then} \texttt{$d_{i,j,k}=d_{new}$}
		\EndIf
		\EndFor
		\EndFor
		\EndFor
		\EndFor
	\end{algorithmic}
\end{algorithm}

\clearpage

%%%%%%%%%%%%%%%%%%%%%%%%%%%%%%%%%%%%%%%%%%%%%%%%%%%%%%%%%%%%%%Vector distance transform%%%%%%%%%%%%%%%%%%%%%%%%%%%%%%%%%%%%%%%%%%%%%%%%%%%%%%%%%%%%%%%

\subsection{Vector distance transform}

For the implementation of the vector distance transform (VDT) \cite{Danielsson} algorithm, we follow the implementation used in \cite{Smisek_dissertation} and extend it to 3D calculations. Comparing the pseudo-code of this method, Alg. \ref{VDT_pseudo_code} with Alg. \ref{FSM_pseudo_code}, we can immediately see that the algorithm also uses Gauss-Seidel iterations alternating the sweeping ordering \eqref{sweeps}. This shows that the information propagates in the same manner, and we can use the same update rules for the values of $d\left ( x \right)$. 
\par The main difference between VDT and FSM lies in the method of how the values of $d\left ( x \right)$ are calculated at the not fixed grid points. While FSM calculates new distance values from the values of neighboring grid points, VDT only checks the source of the neighbors to calculate the smallest possible exact Euclidean distance \eqref{eq:pythagoras} at the current grid point. For this reason, we need to keep track of the sources, and every time we calculate a smaller distance value we update this information. This method yields $O\left(N\right)$ complexity as well.

\begin{algorithm}
	\caption{Vector distance transform}
	\label{VDT_pseudo_code}
	\begin{algorithmic}[1]
		\renewcommand{\algorithmicrequire}{\textbf{Input:}}
		\Require \texttt{From Alg.\ref{Initialization_pc}: 3D grid, $d_{i,j,k}$, $c_{i,j,k}$, $f_{i,j,k}$, $s_{i,j,k}$}
		%\renewcommand{\algorithmicrequire}{\textbf{Definition:}}
		%\mbox{\texttt{~~~~~~$8$ alternating orderings:}}\linebreak
		%\Require \texttt{Eight alternating orderings \eqref{sweeps}. Set of 6 neighbors $P^{1}$ \eqref{set:6_neighbors}.}
		\For{$\left(l=0;l<8;l=l+1\right)$}
		\For{$\left(i=i_{sweep}\left[l,0\right];i \leq i_{sweep}\left[l,1\right];i=i+i_{sweep}\left[l,2\right]\right)$}
		\For{$\left(j=j_{sweep}\left[l,0\right];j \leq j_{sweep}\left[l,1\right];j=j+j_{sweep}\left[l,2\right]\right)$}
		\For{$\left(k=k_{sweep}\left[l,0\right];k \leq k_{sweep}\left[l,1\right];k=k+k_{sweep}\left[l,2\right]\right)$}
		\If{$f_{i,j,k}$ \texttt{is not true} }
		\For{\texttt{all $\left \{ \left( i + r, j + s, k + t \right) ; \left(r,s,t\right) \in P^{1}\right \}$ not out of bound}}
		\If{\texttt{$s_{i+r,j+s,k+t}$ is known}}
		\State $d_{new} = d\left(s_{i+r,j+s,k+t},c_{i,j,k}\right)$	 \Comment{Calculated by \eqref{eq:pythagoras}.}		
		\If{$d_{new}<d_{i,j,k}$}
		\State $d_{i,j,k}=d_{new}$
		\State $s_{i,j,k}=s_{i+r,j+s,k+t}$
		\EndIf
		\EndIf
		\EndFor
		\EndIf
		\EndFor
		\EndFor
		\EndFor
		\EndFor
	\end{algorithmic}
\end{algorithm}

%%%%%%%%%%%%%%%%%%%%%%%%%%%%%%%%%%%%%%%%%%%%%%%%%%%%%%%%%%%%%%Fast marching method%%%%%%%%%%%%%%%%%%%%%%%%%%%%%%%%%%%%%%%%%%%%%%%%%%%%%%%%%%%%%%%

\subsection{Fast marching method}

Similarly, as the FSM algorithm, the fast marching method (FMM) \cite{Sethian_FMM} gives results based on the solution of the Eikonal equation. \textcolor{black}{While FSM tests the possible solutions of the alternatives of \eqref{quadratic2} by going through them in the right order, FMM sets up the solution immediately according to which coefficients are already calculated.} In the construction of this algorithm, one-way propagation of information is utilized, secured by the upwind difference structure of discretization. To properly monitor this propagation the visiting of grid points is tracked throughout the execution of the algorithm. The solution is built outward from the smallest values, which, as seen in the initialization phase in Section \ref{Initialization_phase}, are at the grid points nearest to the points in the cloud. These elements are gathered in a \textit{min-priority-heap} and marked as 'to be visited', while all others are marked 'unvisited'. In Alg. \ref{FMM_pseudo_code} we can see how the heap is used. While the solution from the initialized grid points is marched forward the values from the heap are finalized, marked as 'visited', and new points are brought into this set. FMM works, because we always select the grid point with the smallest value from the heap to calculate the values of the neighboring elements, thus 'unvisited' grid points will not have any effect on the solution.
\par The complexity of the FMM algorithm is of order $O\left(N\,log_{2}N\right)$, because we visit every grid point once and the operations of the \textit{min-priority-heap} have a complexity of $O\left(log_{2}N\right)$.

\begin{algorithm}
	\caption{Fast marching method}
	\label{FMM_pseudo_code}
	\begin{algorithmic}[1]
		\renewcommand{\algorithmicrequire}{\textbf{Input:}}
		\Require \texttt{From Alg.\ref{Initialization_pc}: 3D grid, $d_{i,j,k}$, $f_{i,j,k}$ }
		%\renewcommand{\algorithmicrequire}{\textbf{Definition:}}
		%\Require \texttt{Set of 6 neighbors $P^{1}$ \eqref{set:6_neighbors}.}
		\renewcommand{\algorithmicrequire}{\textbf{Declaration:}}
		\Require \texttt{$v_{i,j,k}$ will hold the visiting values of grid points} \linebreak
		\texttt{'unvisited'=0, 'to be visited'=1, 'visited'=2}
		\Require \texttt{$heap$ container will be a \textit{min-priority-heap}}
		\renewcommand{\algorithmicrequire}{\textbf{Initialization:}}
		\Require \texttt{$\forall f_{i,j,k} = true: \left\lbrace v_{i,j,k} = 1; \; heap.InsertNode(d_{i,j,k}) \right\rbrace $ else: $v_{i,j,k} = 0$}
		\While{\texttt{$heap$ is not empty}}
		\State $\left(i,j,k\right) = heap.GetRoot()$ \Comment{Obtain $(i,j,k)$ with minimum $d$ and delete from $heap$.}
		\For{\texttt{all $\left \{ \left( i + r, j + s, k + t \right) ; \left(r,s,t\right) \in P^{1}\right \}$, not out of bound}}
		\If{\texttt{ ($f_{i + r, j + s, k + t}$ \texttt{is $false$}) and ($v_{i + r, j + s, k + t} = 0 $ or $v_{i + r, j + s, k + t} = 1 $)} }
		\State $x = min\left(d_{i+r+1,j+s,k+t}, d_{i+r-1,j+s,k+t}\right)$ 
		\State $y = min\left(d_{i+r,j+s+1,k+t}, d_{i+r,j+s-1,k+t}\right)$
		\State $z = min\left(d_{i+r,j+s,k+t+1}, d_{i+r,j+s,k+t-1}\right)$ \Comment{Use $+\infty$ if $\left(i+r,j+s,k+t\right)$ is}
		\State $a = b = c = 0$ ~~~~~~~~~~~~~~~~~~~~~~~~~~~~~~~~~~~~~~~~~out of bounds..
		\State \textbf{if} $x \neq +\infty$ \textbf{then} $a = a + 1; \:\: b = b + x; \:\: c = c + x^{2}$
		\State \textbf{if} $y \neq +\infty$ \textbf{then} $a = a + 1; \:\: b = b + y; \:\: c = c + y^{2}$
		\State \textbf{if} $z \neq +\infty$ \textbf{then} $a = a + 1; \:\: b = b + z; \:\: c = c + z^{2}$
		\State $a = a*\left(1/h^{2}\right)$
		\State $b = \left(-2\right)*b*\left(1/h^{2}\right)$
		\State $c = c*\left(1/h^{2}\right) - 1.0$
		\State $d_{new} = \frac{-b + \sqrt{b^{2} - 4 * a*c}}{2*a}$
		\If{$d_{new}<d_{i+r,j+s,k+t}$}
			\State $d_{i+r,j+s,k+t}=d_{new}$
			%\State \textbf{if} \texttt{$v_{i+r,j+s,k+t} = 0 $} \textbf{then} \texttt{$\left\lbrace heap.InsertNode(d_{i+r,j+s,k+t}); \; v_{i+r,j+s,k+t} = 1 \right\rbrace$}
			%\State \textbf{else} $heap.DecreaseKey(\left(i+r,j+s,k+t\right),d_{new})$
			\If{\texttt{$v_{i+r,j+s,k+t} = 0 $}}
				\State $heap.InsertNode(d_{i+r,j+s,k+t})$
				\State	$v_{i+r,j+s,k+t} = 1$
			\Else
				\State $heap.DecreaseKey(\left(i+r,j+s,k+t\right),d_{new})$
			\EndIf
		\EndIf
		\EndIf
		\EndFor
		\State $v_{i,j,k} = 2$
		\EndWhile
	\end{algorithmic}
\end{algorithm}

%%%%%%%%%%%%%%%%%%%%%%%%%%%%%%%%%%%%%%%%%%%%%%%%%%%%%%%%%%%Dijkstra-Pythagoras method%%%%%%%%%%%%%%%%%%%%%%%%%%%%%%%%%%%%%%%%%%%%%%%%%%%%%%%%%%%%

\subsection{Dijkstra-Pythagoras method}

The Dijkstra-Pythagoras (DP) method was introduced in \cite{Smisek_dissertation}. In \cite{Smisek_dissertation}, a gap was detected for a wave-front type method, like FMM, which would yield results with the exact Euclidean distance. Thus the DP method was created. DP algorithm uses visiting rules and a \textit{min-priority-heap} as described in the FMM algorithm but utilizes the source tracking for distance calculation as in the VDT method. In \cite{Smisek_dissertation} the pseudo-code of the method was outlined in a 2D pixel grid with pixel edge size $1$. \textcolor{black}{We extend it to the 3D voxel grid and introduce a substantial modification. In the initial proposal, the algorithm analyzes all neighbors of grid points.  We changed this to include only the closest ones, which in 3D are the voxels from the set  $P^{1}$  \eqref{set:6_neighbors}. We found that with this modification the method becomes much faster and its precision stays approximately the same. In Alg. \ref{DP_pseudo_code} we show the detailed pseudo-code with our changes.} 
\par The logic of the method is based on a two-fold relaxation of $d\left ( x \right)$ values. As in FMM, every cycle of the algorithm starts with the grid point of the smallest $d$ value popped from a \textit{min-priority-heap}. The distance value of this point is checked to the sources of all its ‘visited’ neighbors. From all the $6$ possibilities the value is adjusted to the minimum before it is marked as ‘visited’ as well. Its source is selected accordingly. Then this method attempts to relax the ‘unvisited’ and ‘to be visited’ neighbors in a Dijkstra way. The distances for these grid points are updated according to the Pythagoras rule if the new value is smaller than the value already stored. Their sources are set to the source of the grid point by which they were updated. The neighbors which are ‘unvisited’ will be added to the heap. The algorithm runs till the heap is empty.
\par Similarly to FMM the complexity of this method is $O\left(N\,log_{2}N\right)$.

\clearpage

\begin{algorithm}[]
	\caption{Dijkstra-Pythagoras method}
	\label{DP_pseudo_code}
	\begin{algorithmic}[1]
		\renewcommand{\algorithmicrequire}{\textbf{Input:}}
		\Require \texttt{From Alg.\ref{Initialization_pc}: 3D grid, $d_{i,j,k}$, $c_{i,j,k}$, $f_{i,j,k}$, $s_{i,j,k}$}
		%\renewcommand{\algorithmicrequire}{\textbf{Definition:}}
		%\Require \textcolor{red}{$P^{1}=\left\{\left(r,s,t\right);r,s,t\in\left\{-1,0,1\right\};\left|r\right|+\left|s\right|+\left|t\right|=1\right\}$}
		\renewcommand{\algorithmicrequire}{\textbf{Declaration:}}
		\Require \texttt{$v_{i,j,k}$ will hold the visiting values of grid points} \linebreak
		\texttt{'unvisited'=0, 'to be visited'=1, 'visited'=2}
		\Require \texttt{$heap$ container will be a \textit{min-priority-heap}}
		\renewcommand{\algorithmicrequire}{\textbf{Initialization:}}
		\Require \texttt{$\forall f_{i,j,k} = true: \left\lbrace v_{i,j,k} = 1; \; heap.InsertNode(d_{i,j,k}) \right\rbrace $ else: $v_{i,j,k} = 0$}
		\While{\texttt{$heap$ is not empty}}
			\State $\left(i,j,k\right) = heap.GetRoot()$ \Comment{Obtain $(i,j,k)$ with minimum $d$ and delete from $heap$.}
			\For{\texttt{all $\left \{ \left( i + r, j + s, k + t \right) ; \left(r,s,t\right) \in P^{1}\right \}$, not out of bound}}
				\If{\texttt{$f_{i + r, j + s, k + t}$ \texttt{is $false$} and $v_{i+r,j+s,k+t} = 2$}}
					\State $d_{new} = d\left(s_{i+r,j+s,k+t},c_{i,j,k}\right)$	 \Comment{Calculated by \eqref{eq:pythagoras}.}	
					\If{$d_{new}<d_{i,j,k}$}
						\State $d_{i,j,k}=d_{new}$
						\State $s_{i,j,k}=s_{i+r,j+s,k+t}$
					\EndIf
				\EndIf
			\EndFor
			\State $v_{i,j,k} = 2$
			\For{\texttt{all $\left \{ \left( i + r, j + s, k + t \right) ; \left(r,s,t\right) \in P^{1}\right \}$, not out of bound}}
				\If{\texttt{($f_{i + r, j + s, k + t}$ \texttt{is $false$}) and ($v_{i+r,j+s,k+t} = 0 $ or $v_{i+r,j+s,k+t} = 1 $)}}
					\State $d_{new} = d_{i,j,k} + h$
					\If{$d_{new}<d_{i+r,j+s,k+t}$}
						\State $d_{i+r,j+s,k+t}=d_{new}$
						\State $s_{i+r,j+s,k+t}=s_{i,j,k}$
						%\State \textbf{if} \texttt{$v_{i+r,j+s,k+t} = 0 $} \textbf{then} \texttt{$\left\lbrace heap.InsertNode(d_{i+r,j+s,k+t}); \; v_{i+r,j+s,k+t} = 1 \right\rbrace$}
						%\State \textbf{else} $heap.DecreaseKey(\left(i+r,j+s,k+t\right),d_{new})$
						\If{\texttt{$v_{i+r,j+s,k+t} = 0 $}}
							\State $heap.InsertNode(d_{i+r,j+s,k+t})$
							\State	$v_{i+r,j+s,k+t} = 1$
						\Else
							\State $heap.DecreaseKey(\left(i+r,j+s,k+t\right),d_{new})$
						\EndIf
					\EndIf
				\EndIf
			\EndFor
		\EndWhile
	\end{algorithmic}
\end{algorithm}

%%%%%%%%%%%%%%%%%%%%%%%%%%%%%%%%%%%%%%%%%%%%%%%%%%%%%%%%%%%Comparing_methods%%%%%%%%%%%%%%%%%%%%%%%%%%%%%%%%%%%%%%%%%%%%%%%%%%%%%%%%%%%%

\section{Numerical experiments - methods comparison} \label{Numerical_experiments}

In this section, we compare the efficiency of the described algorithms and show they can be used for computing the distance function to objects represented by a 3D point cloud and triangulated surface.

%%%%%%%%%%%%%%%%%%%%%%%%%%%%%%%%%%%%%%%%%%%%%%%%%%%%%%%%%%%Comparing methods%%%%%%%%%%%%%%%%%%%%%%%%%%%%%%%%%%%%%%%%%%%%%%%%%%%%%%%%%%%%

\subsection{Comparing methods}\label{sec:comparing_method}

\textcolor{black}{For the first experiment, we will work with a cube with an edge size of 1.0. an its vertex with minimum coordinates at $\left(0.0,0.0,0.0\right)$. We construct around it a rectangular computational domain which is $0.4$ times larger from the Cube in every direction. In this experiment, we discretize the computational domain in a way that some of the grid points will always lie on the surface of the Cube. Thus, we can set the distance function at these points to $0$ during initialization.}
\par
\begin{figure}[htp]
	\centering
	\begin{minipage}{0.49\linewidth}
		\centering
		\includegraphics[width=\linewidth]{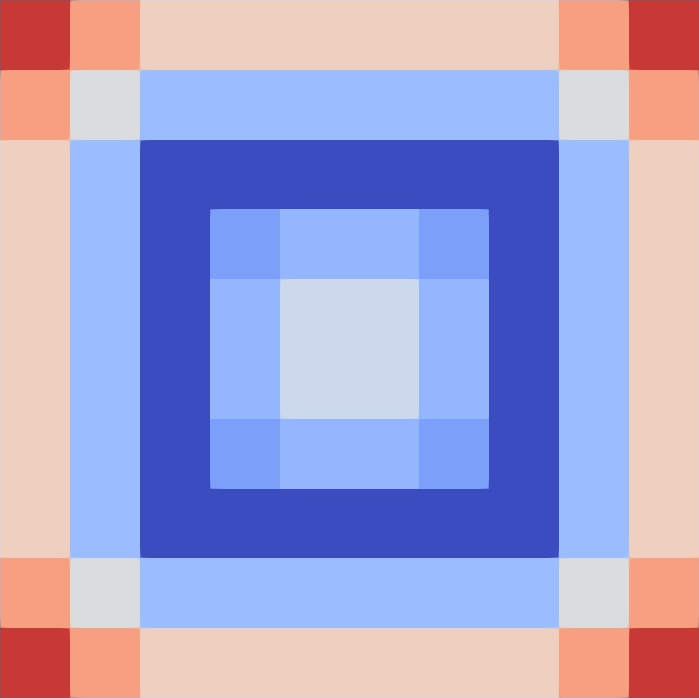}
	\end{minipage}
	\vspace{0.04cm}
	\begin{minipage}{0.49\linewidth}
		\centering
		\includegraphics[width=\linewidth]{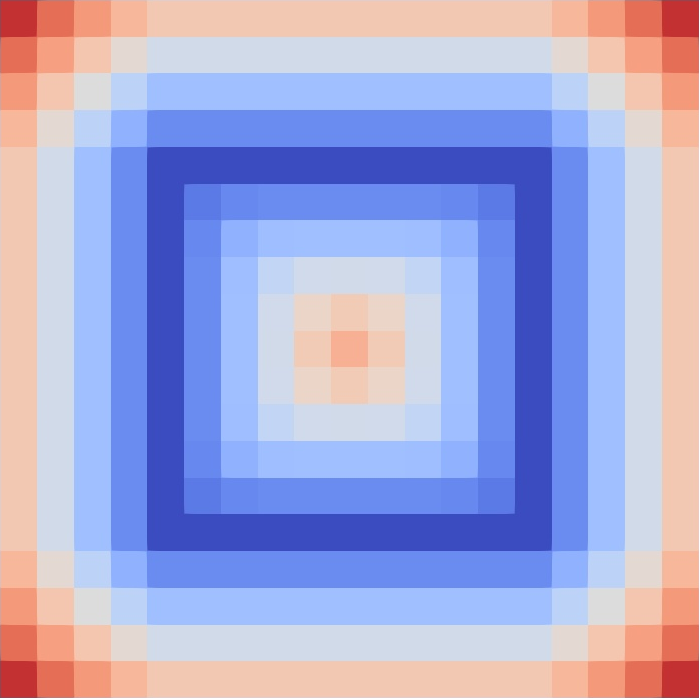}
	\end{minipage}
	\vspace{0.04cm}
	
	\begin{minipage}{0.49\linewidth}
		\centering
		\includegraphics[width=\linewidth]{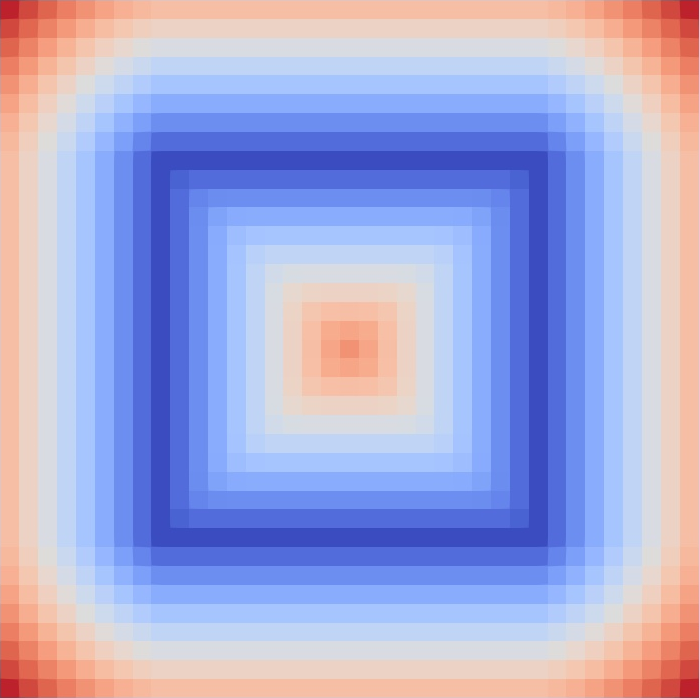}
	\end{minipage}
	\vspace{0.04cm}
	\begin{minipage}{0.49\linewidth}
		\centering
		\includegraphics[width=\linewidth]{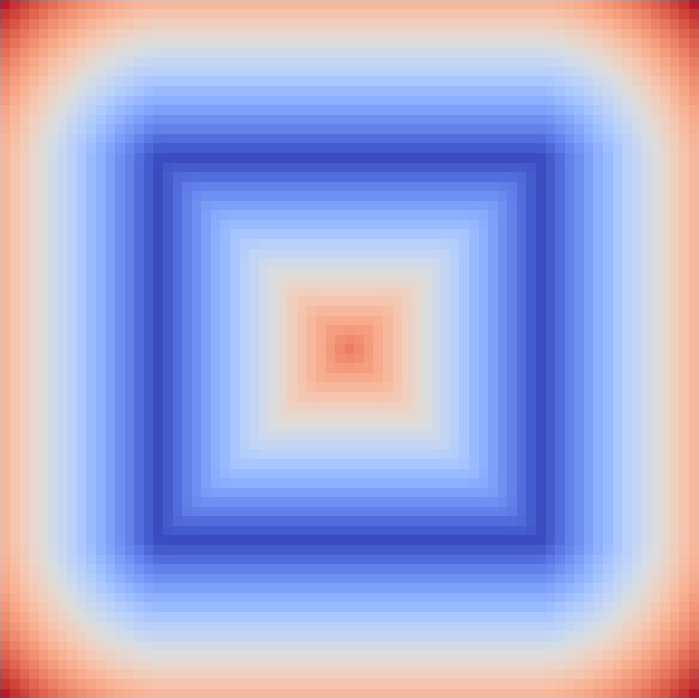}
	\end{minipage}
	\vspace{0.04cm}
	
	\caption{Distance function visualization for the Cube experiment. We visualize the section in a constant $z$ plane for voxel edge sizes $0.2$, $0.1$, $0.05$, $0.025$. Values go from the highest dark red to the lowest dark blue. Results were calculated by FSM.}
	\label{fig:cube_distance}
\end{figure}

 With this setup, we computed the distance function for the Cube with the four algorithms on the computational domain discretized to a grid by voxels with different edge sizes, namely $0.2$, $0.1$, $0.05$, $0.025$, $0.0125$, $0.00625$, $0.003125$. We demonstrate how the distance function looks like on these grids in Figure \ref{fig:cube_distance} calculated by the FSM algorithm. 

\clearpage

\par
\begin{figure}[htp]
	\centering
	
	\begin{minipage}{0.24\linewidth}
		\centering
		\includegraphics[width=\linewidth]{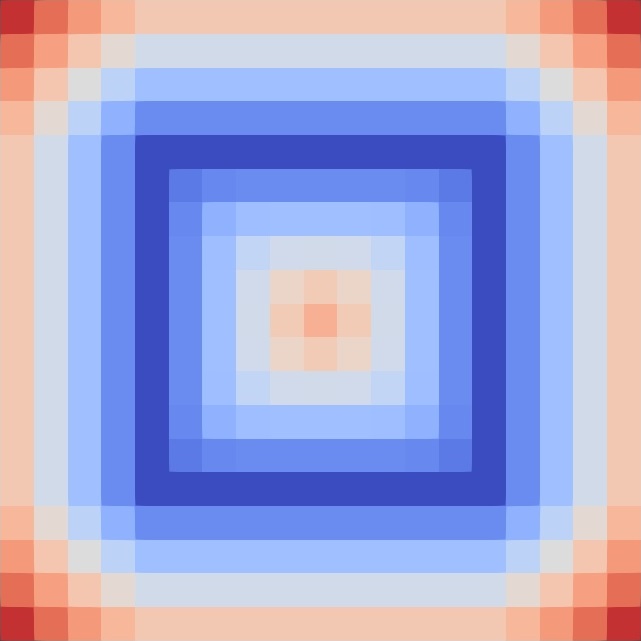}
	\end{minipage}
	\vspace{0.035cm}
	\begin{minipage}{0.24\linewidth}
		\centering
		\includegraphics[width=\linewidth]{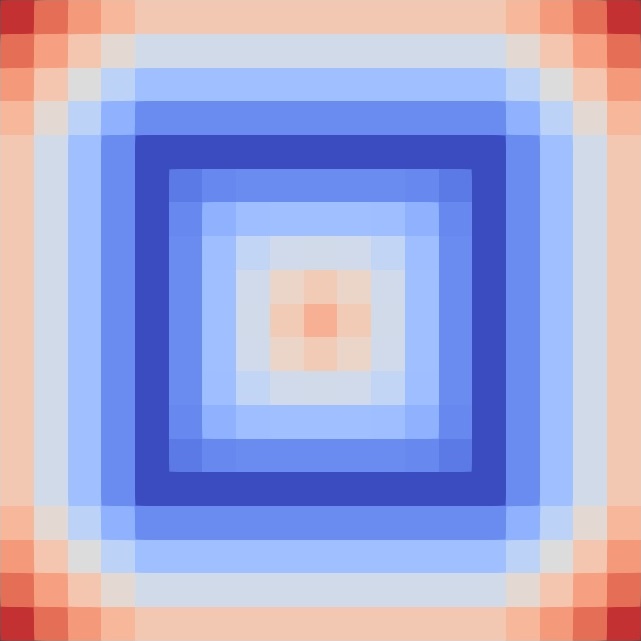}
	\end{minipage}
	\vspace{0.035cm}
	\begin{minipage}{0.24\linewidth}
		\centering
		\includegraphics[width=\linewidth]{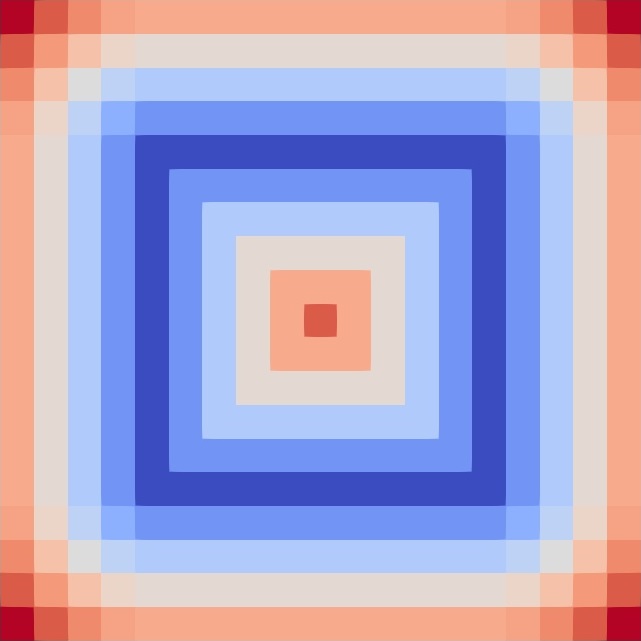}
	\end{minipage}
	\vspace{0.035cm}
	\begin{minipage}{0.24\linewidth}
		\centering
		\includegraphics[width=\linewidth]{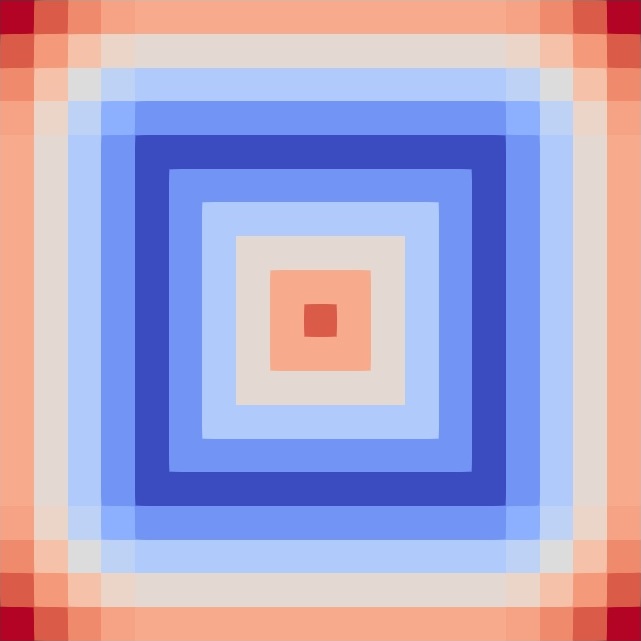}
	\end{minipage}
	\vspace{0.035cm}
	
	\begin{minipage}{0.24\linewidth}
		\centering
		\includegraphics[width=\linewidth]{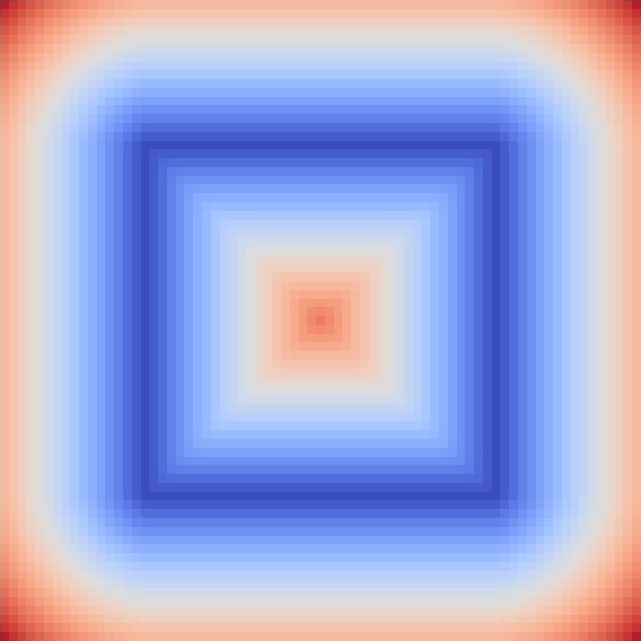}
	\end{minipage}
	%\vspace{0.035cm}
	\begin{minipage}{0.24\linewidth}
		\centering
		\includegraphics[width=\linewidth]{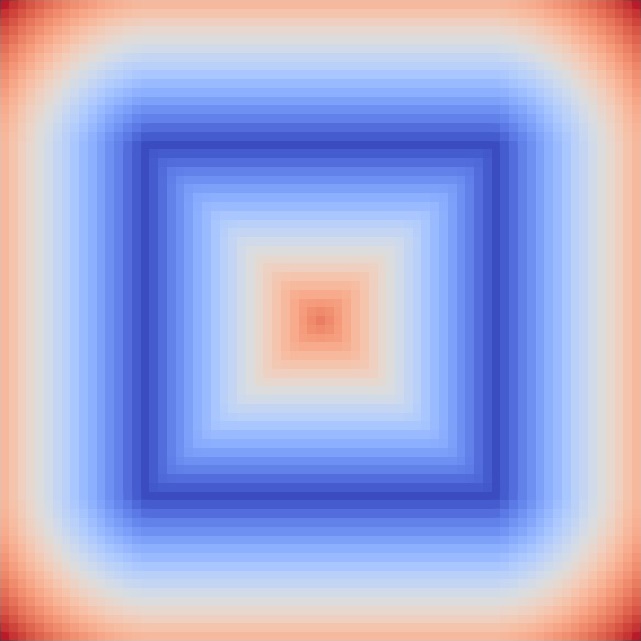}
	\end{minipage}
	%\vspace{0.035cm}
	\begin{minipage}{0.24\linewidth}
		\centering
		\includegraphics[width=\linewidth]{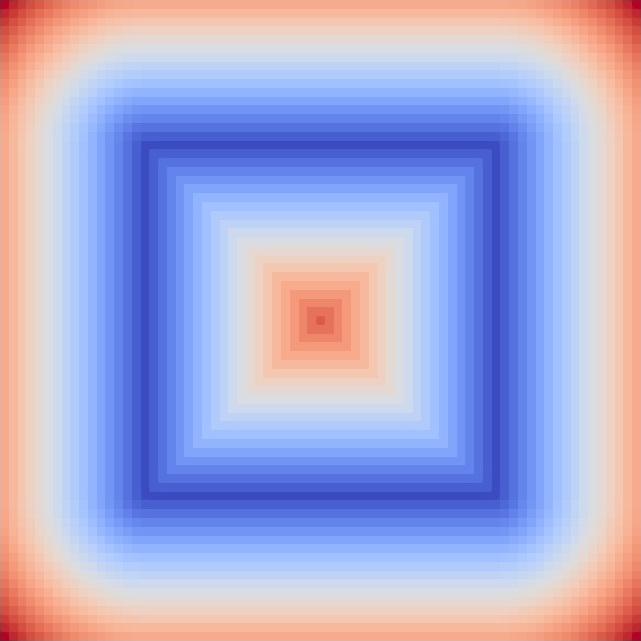}
	\end{minipage}
	%\vspace{0.035cm}
	\begin{minipage}{0.24\linewidth}
		\centering
		\includegraphics[width=\linewidth]{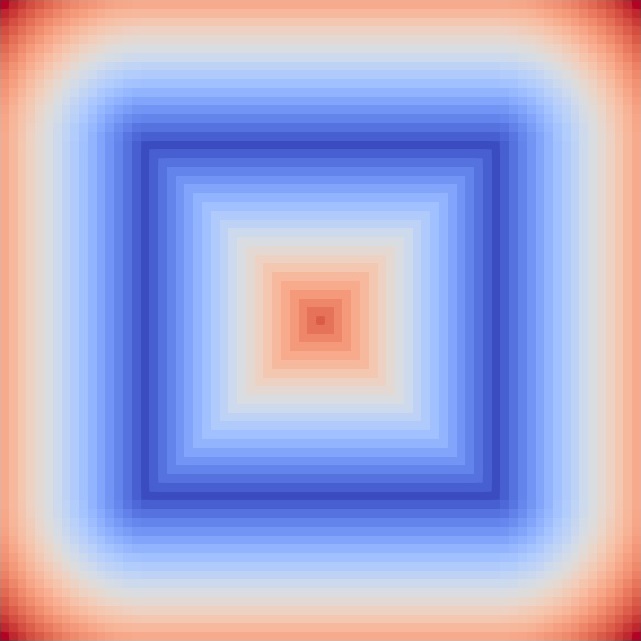}
	\end{minipage}
	%\vspace{0.035cm}
	
	\caption{Visualization of results for distance function calculation in a constant $z$ plane. In the first row, we see visualization for voxel edge size $0.1$, in the second row for voxel edge size $0.25$. In the first column, we see the result for the FSM algorithm, in second for FMM, in third for VDT and in the fourth for DP.}
	\label{fig:Cube_DisComp}
\end{figure}

 To compare the accuracy of the algorithms, we calculated the mean squared difference from the exact solution for all grids. If we denote the exact solution as $\bar{d}_{i,j,k}$ at $x_{i,j,k} \in \Omega$ the mean squared difference can be calculated as
\begin{equation} \label{eq:mean_squared_diff}
	\left [ \sum_{i = 0}^{N_{1}-1} \sum_{j = 0}^{N_{2}-1} \sum_{k = 0}^{N_{3}-1} \left (\bar{d}_{i,j,k} - d_{i,j,k}\right )^{2} \right ] / \left (N_{1}*N_{2}*N_{3}\right )
\end{equation}
We are listing these results in Table \ref{tab:cube_error}. In the first column, we list the number of grid points in $x$, $y$, $z$ directions of our computational grid. In the second column, we see the length of the voxel edges. In the following columns, we see the mean squared difference for FSM, VDT, FMM, and DP methods. We can see that the results for the VDT and DP methods are basically $0$, as we would have expected after stating the fact that they yield Euclidean distance results. The results of FSM and FMM are less accurate. We compare these results also visually in Figure \ref {fig:Cube_DisComp} for computational grids with voxel edge size $0.1$ and $0.025$. We can see that the results for the pairs of FSM, FMM, and VDT, DP in this experiment are visually identical.
\par Besides the accuracy, for this experiment, we also measured the CPU time in seconds which was needed to calculate the distance function with the different methods, reported in Table \ref{tab:cube_cputime}. Here again, we list the parameters of our grid first. In the third column, we list the CPU time for the initialization phase of the algorithms. The initialization is the same for all four methods. Because of the simplicity of the experiment, this takes just a few seconds even for the finest grid. Comparing the results we see that concerning CPU time the FSM algorithm outperforms all other methods.

\begin{table}[htp]
	\centering
	\begin{tabular}{|c|c|c|c|c|c|}
		\hline
		\begin{tabular}[c]{@{}c@{}}Number of \\  grid points\end{tabular} & \begin{tabular}[c]{@{}c@{}}Voxel\\ edge size\end{tabular} & FSM        & VDT        & FMM        & DP                              \\ \hline
		$10^{3}$	& 0.2		& 2.5692e-03 	& 1.4791e-34 	& 2.5692e-03 	& 4.227801e-33	\\ \hline
		$19^{3}$    & 0.1   	& 9.7901e-04 	& 1.8869e-34 	& 9.7902e-04 	& 5.011068e-33 	\\ \hline
		$37^{3}$    & 0.05  	& 3.7697e-04 	& 1.0579e-34 	& 3.7697e-04 	& 1.151981e-32  \\ \hline
		$73^{3}$    & 0.025 	& 1.4352e-04 	& 5.0275e-35 	& 1.4352e-04 	& 8.414407e-33  \\ \hline
		$145^{3}$   & 0.0125    & 5.3092e-05 	& 2.5195e-35 	& 5.3092e-05 	& 9.454389e-33  \\ \hline
		$289^{3}$   & 0.00625   & 1.8949e-05 	& 1.3057e-35 	& 1.8949e-05 	& 3.427338e-32  \\ \hline
		$577^{3}$   & 0.003125  & 6.5244e-06 	& 6.5770e-36 	& 6.5244e-06 	& 1.404635e-31 	\\ \hline
	\end{tabular}
	\caption{Mean squared difference comparison for distance function calculation methods tested on the Cube experiment.}
	\label{tab:cube_error}
\end{table}

\clearpage

\begin{table}[htp]
	\centering
	\begin{tabular}{|c|c|c|c|c|c|c|}
		\hline
		\begin{tabular}[c]{@{}c@{}}Number of \\ grid points\end{tabular} & \begin{tabular}[c]{@{}c@{}}Voxel \\ edge size\end{tabular} & Initialization & FSM    & VDT    & FMM    & DP      \\ \hline
		$10^{3}$	& 0.2		& 0			& 0.001  & 0.002  & 0.001  	& 0.001   \\ \hline
		$19^{3}$	& 0.1  		& 0			& 0.002  & 0.017  & 0.002  	& 0.002   \\ \hline
		$37^{3}$    & 0.05  	& 0.001     & 0.008  & 0.029  & 0.016  	& 0.015   \\ \hline
		$73^{3}$    & 0.025     & 0.014     & 0.059  & 0.186  & 0.18   	& 0.138   \\ \hline
		$145^{3}$   & 0.0125    & 0.105     & 0.352  & 1.416  & 1.988   & 1.441   \\ \hline
		$289^{3}$   & 0.00625   & 0.816     & 2.881  & 11.352 & 26.109 	& 15.425  \\ \hline
		$577^{3}$   & 0.003125  & 6.375     & 24.442 & 88.836 & 313.009 & 159.762 \\ \hline
	\end{tabular}
	\caption{CPU time comparison for distance function calculation methods tested on the Cube experiment. CPU time was measured in seconds.}
	\label{tab:cube_cputime}
\end{table}

\par
\begin{figure}[htp]
	\centering
	\includegraphics[width=0.49\linewidth]{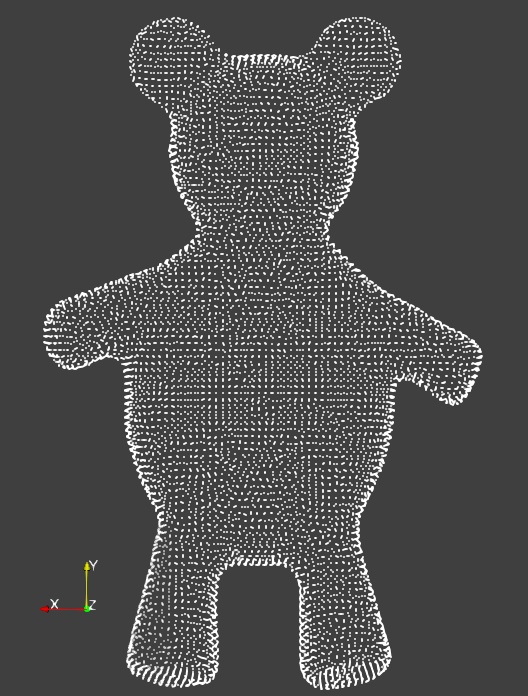}
	\includegraphics[width=0.49\linewidth]{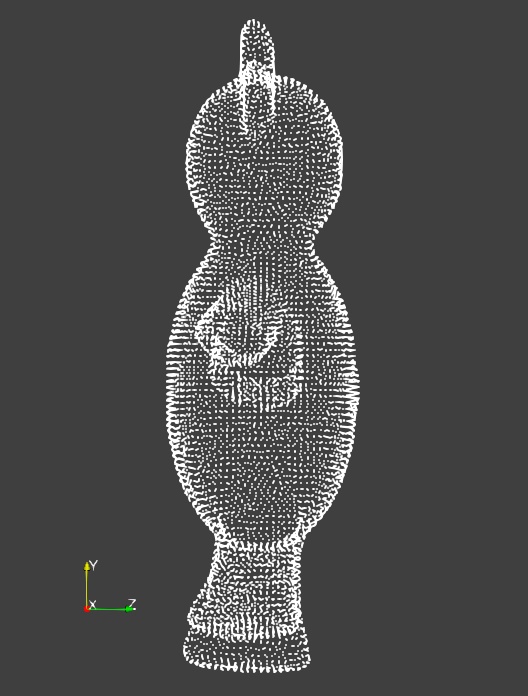}
	\caption{Teddy Bear point cloud data. In the left picture, we can see it from the front in the right picture from the side.}
	\label{fig:bear_pointcloud}
\end{figure}

\begin{figure}[htp]
	\begin{minipage}{0.325\linewidth}
		\centering
		\includegraphics[width=\linewidth]{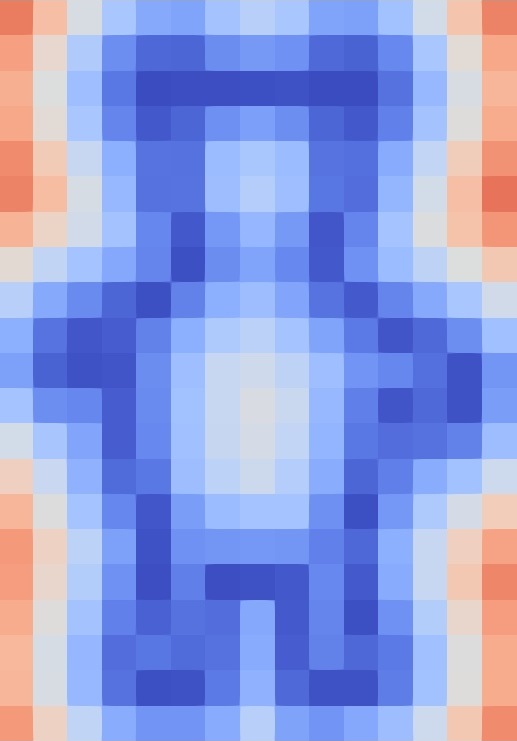}
	\end{minipage}
	\vspace{0.025cm}
	\begin{minipage}{0.325\linewidth}
		\centering
		\includegraphics[width=\linewidth]{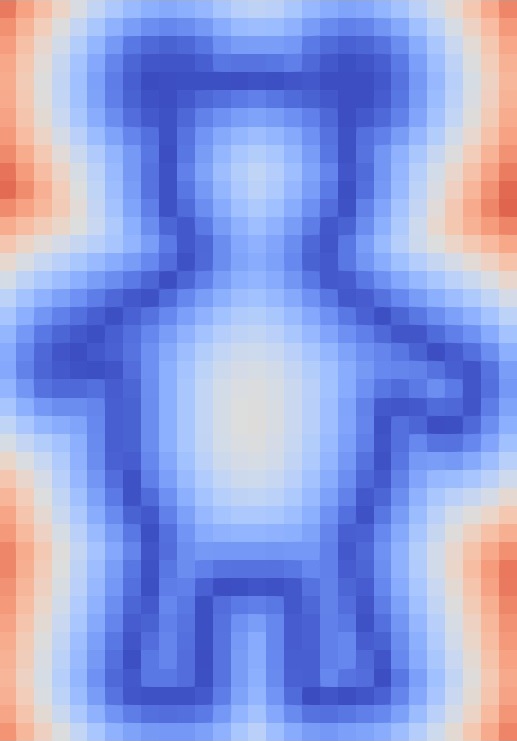}
	\end{minipage}
	\vspace{0.025cm}
	\begin{minipage}{0.325\linewidth}
		\centering
		\includegraphics[width=\linewidth]{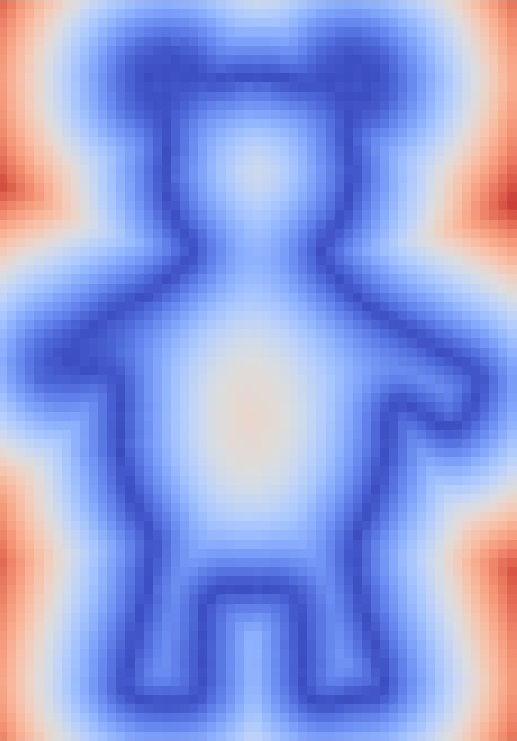}
	\end{minipage}
	\vspace{0.025cm}
	
	\begin{minipage}{0.325\linewidth}
		\centering
		\includegraphics[width=\linewidth]{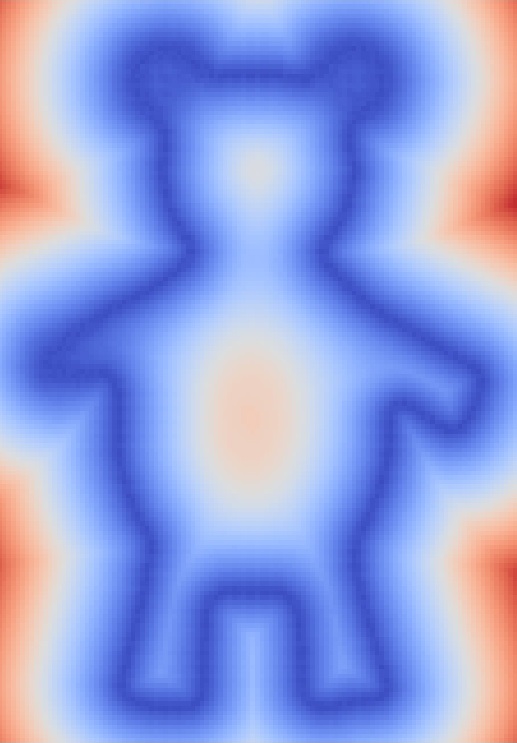}
	\end{minipage}
	\vspace{0.025cm}
	\begin{minipage}{0.325\linewidth}
		\centering
		\includegraphics[width=\linewidth]{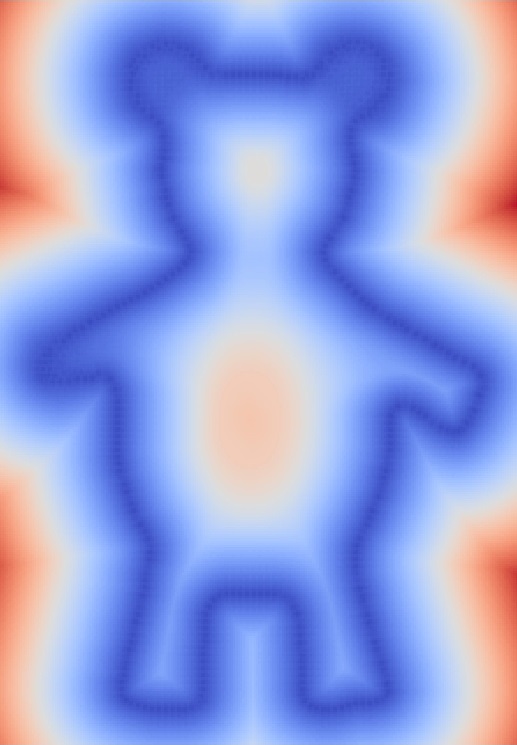}
	\end{minipage}
	\vspace{0.025cm}
	\begin{minipage}{0.325\linewidth}
		\centering
		\includegraphics[width=\linewidth]{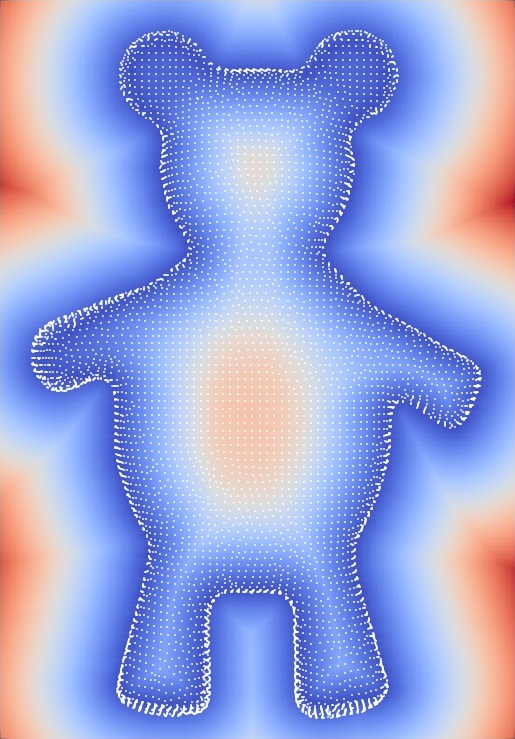}
	\end{minipage}
	\vspace{0.025cm}
	
	\caption{Distance function visualization of the Teddy Bear data set. Visualizing sections in a constant $z$ plane for voxel edge sizes $0.1$, $0.05$, $0.025$, $0.0125$, $0.00625$, $0.003125$. In the last picture, we visualize the distance function with point cloud data. Values go from the highest dark red to the lowest dark blue. Results were calculated by FSM.}
	\label{fig:bear_distance_density}
\end{figure}

 For further comparison of efficiency we choose a data set from \cite{Chen_2009_ABF} which will be used as a point cloud data and as a triangulated surface as well. This data set, seen in Figure \ref{fig:bear_pointcloud}, represents a teddy bear. Similarly, as in the previous experiment, we computed the distance function for the point cloud data with the four algorithms on computational grids with different voxel edge sizes $0.1$, $0.05$, $0.025$, $0.0125$, $0.00625$, $0.003125$, $0.0015625$. Some of the results for distance functions calculated by the FSM algorithm can be seen in Figure \ref{fig:bear_distance_density}. Visually there is no big difference between the results of the four algorithms. 
\par We list the CPU time for calculation in Table \ref{tab:bear_pc_cpu_time}. We added one more information in this table that was not listed in the previous experiment. In the third column, we list the number of fixed grid points produced by the initialization phase of the calculations. We will use this information for the comparison of distance function calculation in the case of the triangulated surface. In this experiment the points from the point cloud data do not coincide with points of the grid, thus the initialization was done by Alg. \ref{Initialization_pc}. The FSM algorithm is the fastest in this case as well.

\begin{table}[htp]
	\centering
	\catcode`\-=12
	\scalebox{0.9}{
		\begin{tabular}{|c|c|c|c|c|c|c|c|}
			\hline
			\begin{tabular}[c]{@{}c@{}}Number of \\  grid points\end{tabular} & \begin{tabular}[c]{@{}c@{}}Voxel\\ edge size\end{tabular} & \begin{tabular}[c]{@{}c@{}}Fixed \\ points\end{tabular} & \begin{tabular}[c]{@{}c@{}}Initial\\ condition\end{tabular} & FSM   & VDT   & FMM   & DP    \\ \hline
			15 x 21 x 9			& 0.1		& 1003		& 0.001		& 0.002 	& 0.002 	& 0.001 	& 0.001   \\ \hline
			29 x 41 x 17		& 0.05  	& 3906  	& 0.002 	& 0.007 	& 0.024 	& 0.005 	& 0.005   \\ \hline
			57 x 81 x 32    	& 0.025 	& 14428 	& 0.007 	& 0.036 	& 0.083 	& 0.059 	& 0.054   \\ \hline
			113 x 161 x 62  	& 0.0125    & 47575 	& 0.042 	& 0.253 	& 0.562 	& 0.658 	& 0.605   \\ \hline
			224 x 321 x 122 	& 0.00625   & 76231 	& 0.293 	& 1.891 	& 4.135 	& 9.079 	& 7.638  \\ \hline
			447 x 640 x 242 	& 0.003125  & 76384 	& 2.235 	& 15.335 	& 33.126 	& 121.185 	& 103.515 \\ \hline
			893 x 1279 x 482 	& 0.0015625 & 76384     & 19.085 	& 118.572 	& 254.895 	& 1439.6	& 1221.13 \\ \hline
	\end{tabular}}
	\caption{CPU time comparison for distance function calculation methods tested on the Teddy Bear point cloud data. CPU time was measured in seconds.}
	\label{tab:bear_pc_cpu_time}
\end{table}

\clearpage

%%%%%%%%%%%%%%%%%%%%%%%%%%%%%%%%%%%%%%%%%%%%%Distance function calculation to triangulated surfaces%%%%%%%%%%%%%%%%%%%%%%%%%%%%%%%%%%%%%%%%%%%%%%%%%

\subsection{Distance function to triangulated surfaces}\label{sec:distance_to_tm}

With small changes, it is possible to easily modify the algorithm for the calculation of the distance function to triangulated surfaces. The most important changes which need to be applied concern the initialization phase. We demonstrate this in the pseudo-code Alg. \ref{Initialization_tm}. In this algorithm, we cycle through all triangles in the triangulated surface. For every triangle, we find the grid points which are lying next to its surface. In these grid points, we calculate the distance from the triangle. For this, we use the method described in \cite{Eberly_PointToTriangle}. Similarly, as with the point cloud data, the values in these points will be fixed, but now as the source of distance computation, we will refer to the triangles. Regarding the algorithms FSM, VDT, FMM, and DP, the only changes will be in the pseudo-code Alg. \ref{VDT_pseudo_code} for VDT on line $8$ and in the pseudo-code Alg. \ref{DP_pseudo_code} for DP on line $5$ where the distance will be calculated between a point and a triangle.

\begin{algorithm}
	\caption{Initialization of distance function to triangulated surface}
	\label{Initialization_tm}
	\begin{algorithmic}[1]
		\renewcommand{\algorithmicrequire}{\textbf{Input:}}
		\Require \texttt{Triangulated surface:} \linebreak
		\texttt{$tr_{l}$ - set of triangles,} \linebreak
		\texttt{$N$ - number of triangles.}
		%\mbox{\texttt{~~~3D grid with voxel edge size $h$ and dimensions:}} \linebreak 
		\Require \texttt{3D grid with voxel edge size $h$ and dimensions $N_{i}, N_{j}, N_{k}$.}
		\renewcommand{\algorithmicrequire}{\textbf{Declaration:}}
		\Require \texttt{Arrays:} \linebreak 
		\texttt{$d_{i,j,k}$ - value of distance function at grid point $\left(i,j,k\right)$, } \linebreak
		\texttt{$c_{i,j,k}$ - $\left(x,y,z\right)$ coordinates of grid point $\left(i,j,k\right)$, } \linebreak
		\texttt{$f_{i,j,k}$ - determines if $d_{i,j,k}$ is fixed at $(i,j,k)$, } \linebreak
		\texttt{$s_{i,j,k}$ - source for $d_{i,j,k}$ calculation at $(i,j,k)$. }
		\State \textbf{Set} \texttt{$d_{i,j,k}$ to $+\infty$, $f_{i,j,k}$ to $false$, $s_{i,j,k}$ to $unknown$ }
		\State \textbf{Calculate:} $c_{i,j,k}$
		\For{$\left(l=0;l<N;l=l+1\right)$}
		\State $gp_{m} = PointsAlongTriangle\left( tr_{l} \right)$ \Comment{$gp_{m}$ is a subset of the computational grid.}
		\State $N_{gp} = NumberOfPointsIn(gp_{m})$
		\For{$\left(m=0;m < N_{gp}; m=m+1\right)$}
		\State $\left(i,j,k\right) = gp_{m}$
		\State $d_{new} = d\left(tr_{l},c_{i,j,k}\right)$	 \Comment{Distance of a point from a triangle.}
		\If{$d_{new}<d_{i,j,k}$}
		\State $d_{i,j,k}=d_{new}$ 
		\State $f_{i,j,k}=true$ 
		\State $s_{i,j,k}=tr_{l}$
		\EndIf
		\EndFor	
		\EndFor
	\end{algorithmic}
\end{algorithm}

\par
\begin{figure}[htp]
	\centering
	\includegraphics[width=0.48\linewidth]{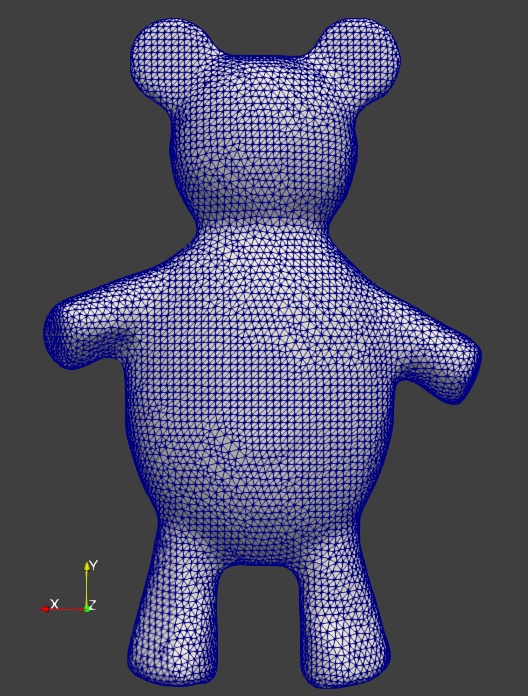}
	\includegraphics[width=0.48\linewidth]{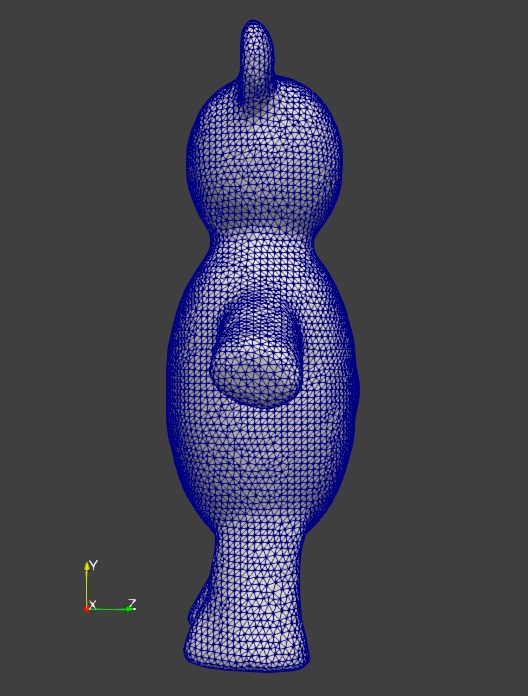}
	\caption{Teddy Bear triangulated surface data. In the left picture, we can see it from the front in the right picture from the side.}
	\label{fig:bear_triangulated}
\end{figure}

\begin{figure}[htp]
	\centering
	\begin{minipage}{0.48\linewidth}
		\centering
		\includegraphics[width=\linewidth]{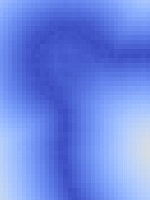}
	\end{minipage}
	\vspace{0.04cm}
	\begin{minipage}{0.48\linewidth}
		\centering
		\includegraphics[width=\linewidth]{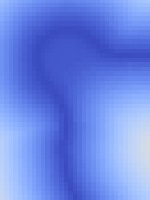}
	\end{minipage}
	\vspace{0.04cm}
	
	\begin{minipage}{0.48\linewidth}
		\centering
		\includegraphics[width=\linewidth]{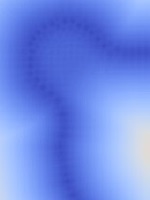}
	\end{minipage}
	\vspace{0.04cm}
	\begin{minipage}{0.48\linewidth}
		\centering
		\includegraphics[width=\linewidth]{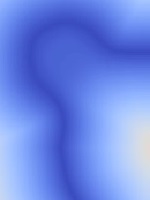}
	\end{minipage}
	\vspace{0.04cm}
	
	\caption{Comparing the results of distance function calculation from point cloud data (first column) and triangular surface (second column). Voxel edge size for results in the first row is $0.0125$, in second row $ 0.003125$. Results were calculated by FSM.}
	\label{fig:bear_pc_vs_tm}
\end{figure}

 To demonstrate the results of these changes we will use again the Teddy Bear data set, but now as a triangulated surface as seen in Figure \ref{fig:bear_triangulated}. Similarly, as for the calculation to the point cloud data, we measured the CPU times and listed them in Table \ref{tab:bear_tm_cpu_time}.
If we compare this to Table \ref{tab:bear_pc_cpu_time} we can see the difference between the calculation of the distance function for point cloud data and a triangulated surface. The number of fixed points is much higher for the triangulated surface. This is because the initialization produces a "contiguous" volume around the triangles for every density of the grid, while around the point cloud data gaps can develop. We can see this also in Figure \ref{fig:bear_pc_vs_tm}. Here we compare the distance function for both point cloud and triangulated surface by the results obtained by the FSM algorithm. (The difference in the visualization of the distance function calculated with the other algorithms is very small thus we provide just the visualization of the FSM algorithm.) The results for the triangulated surface, seen in the right column, are much smoother near the object as for the point cloud data, seen in the left column. While this difference has no real effect on the calculation time of the FSM and FMM algorithms, it drastically increases the time for VDT and DP. This is because the implementation of FSM and FMM is independent of the initial data, but in VDT and DP we work with the source as well and the calculation of the distance between a point and a triangle takes more time than the calculation between two points.
\par
\begin{figure}[htp]
	\centering
	\includegraphics[height=65mm]{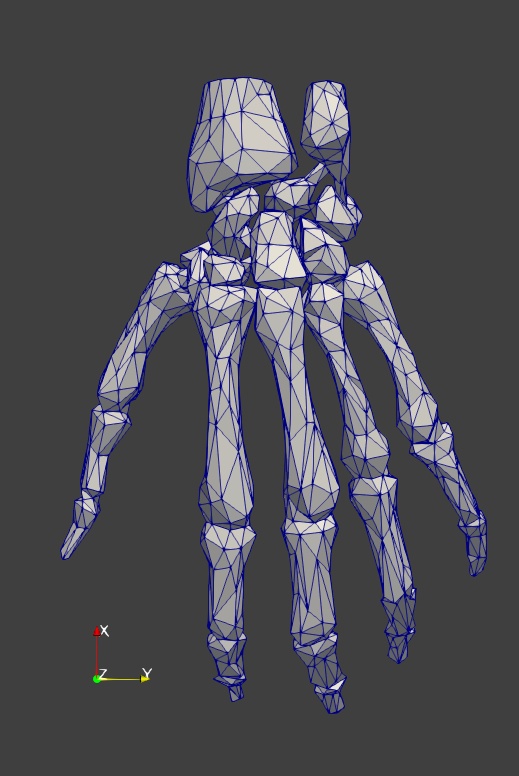}
	\includegraphics[height=65mm]{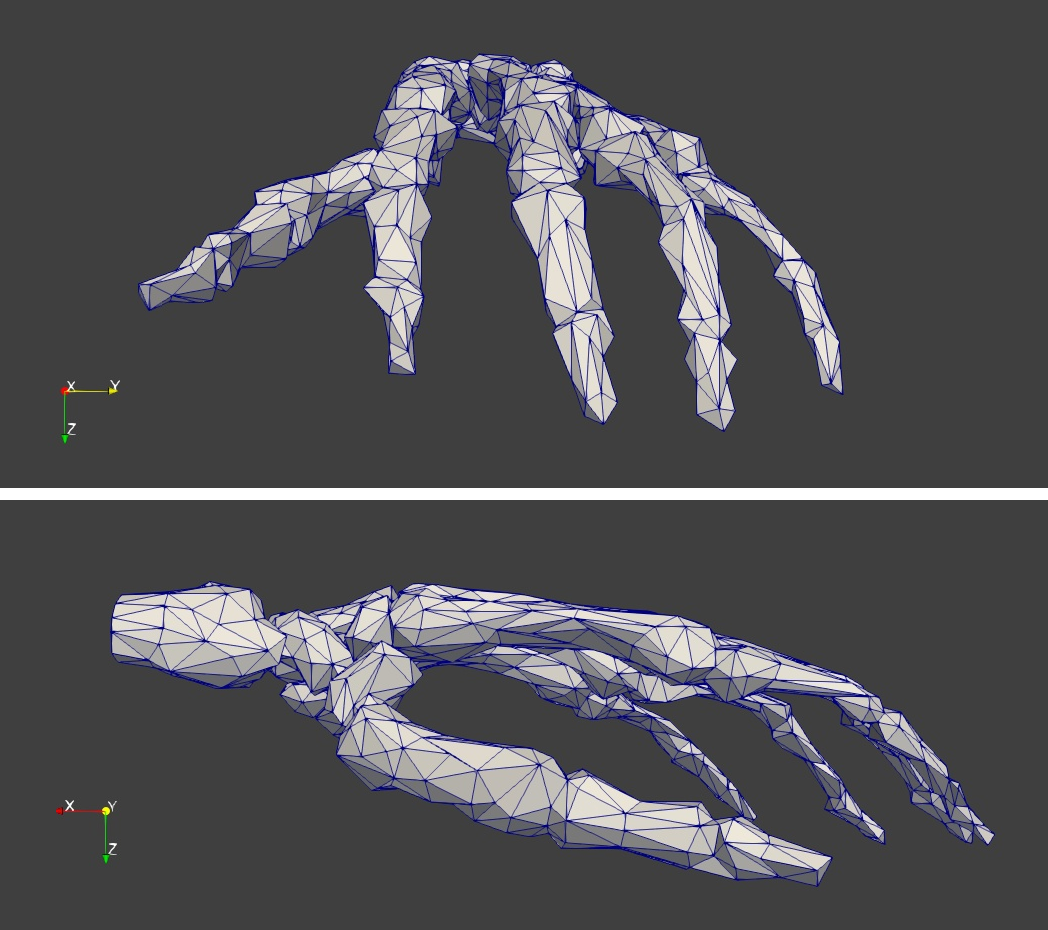}
	\caption{Hand Bones triangulated surface data. In the left picture, we can see it from above, in the right upper picture from the front and in the right bottom picture from the side.}
	\label{fig:hand_skeleton_triangulated}
\end{figure}

\begin{figure}[htp]
	\centering
	\begin{minipage}{\linewidth}
		\centering
		\includegraphics[width=0.49\linewidth]{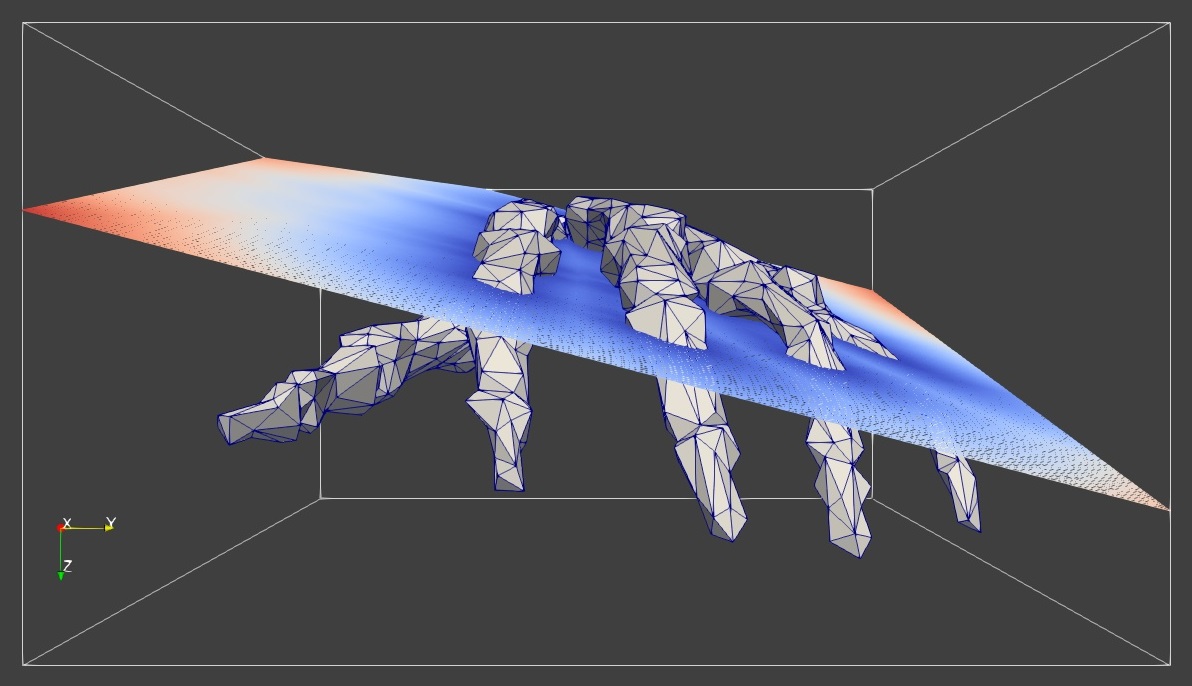}
		\includegraphics[width=0.49\linewidth]{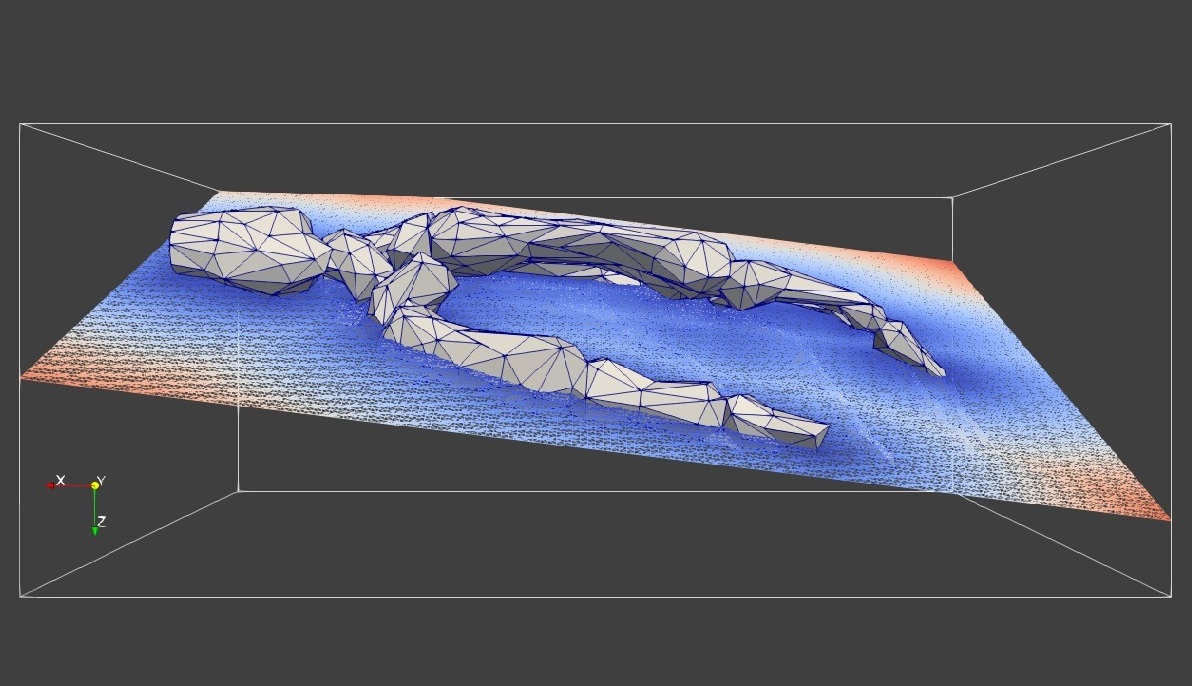}
	\end{minipage}
	\hspace{0.08cm}
	
	\begin{minipage}{\linewidth}
		\centering
		\includegraphics[width=0.49\linewidth]{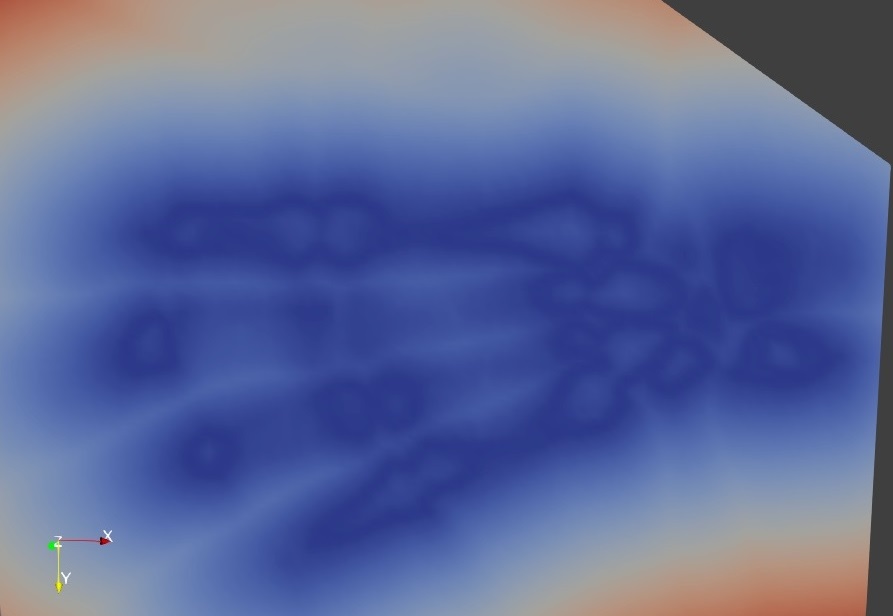}
		\includegraphics[width=0.49\linewidth]{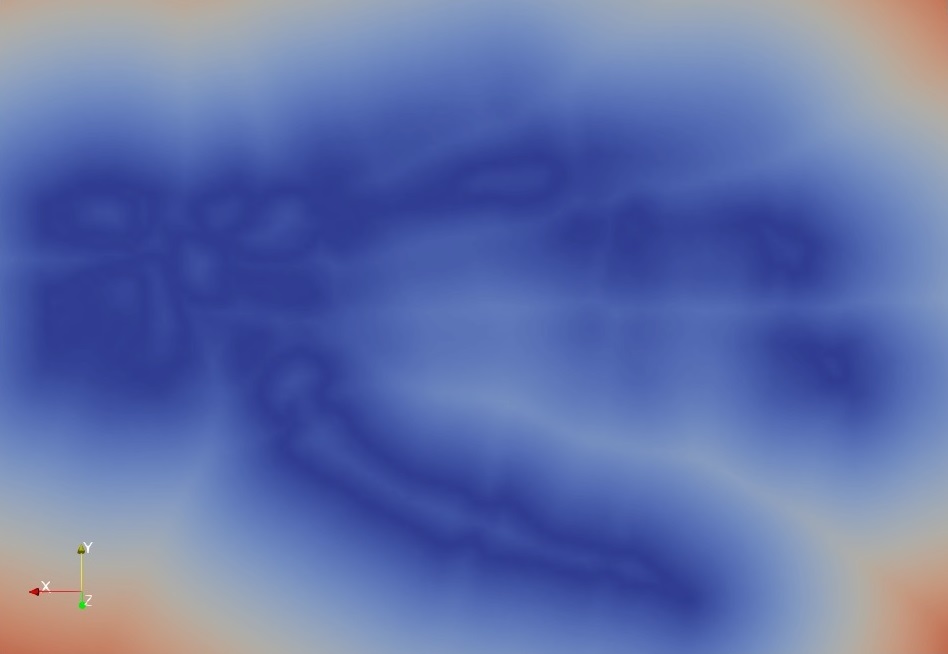}
	\end{minipage}
	\hspace{0.08cm}
	
	\caption{Visualization of slices of distance function calculated to Hand Bones triangulated surface data. In the upper pictures, we can see the location of slices in the 3D computational domain, in the bottom picture the slices.}
	\label{fig:hand_skeleton_slices}
\end{figure}

 \textcolor{black}{To demonstrate a further example of distance function calculation on a triangulated surface we applied the algorithm on an additional data set. We obtained it from \cite{Models_Archive}. In Figure \ref{fig:hand_skeleton_triangulated} we can see the triangulated surface of hand bones. With its many details and small parts, it is a good data set to show the accuracy of the results. We can see these in Figure \ref{fig:hand_skeleton_slices}. Here we choose planes in the computational domain in which we can see the most details.}

\begin{table}[htp]
	\centering
	\catcode`\-=12
	\scalebox{0.9}{
	\begin{tabular}{|c|c|c|c|c|c|c|c|}
	\hline
	\begin{tabular}[c]{@{}c@{}}Number of \\  grid points\end{tabular} & \begin{tabular}[c]{@{}c@{}}Voxel\\ edge size\end{tabular} & \begin{tabular}[c]{@{}c@{}}Fixed \\ points\end{tabular} & \begin{tabular}[c]{@{}c@{}}Initial\\ condition\end{tabular} & FSM   & VDT   & FMM   & DP    \\ \hline
			15 x 21 x 9           & 0.1        & 1283		& 0.039		& 0.001		& 0.007		& 0.001		& 0.001	\\ \hline
			29 x 41 x 17          & 0.05       & 5147		& 0.048		& 0.007		& 0.056		& 0.005		& 0.008	\\ \hline
			57 x 81 x 32          & 0.025      & 20526		& 0.067		& 0.04		& 0.446		& 0.057		& 0.08	\\ \hline
			113 x 161 x 62        & 0.0125     & 82756		& 0.139		& 0.255		& 3.408		& 0.659		& 0.828	\\ \hline
			224 x 321 x 122       & 0.00625    & 331874		& 0.504		& 1.966		& 26.054	& 9.475		& 11.02	\\ \hline
			447 x 640 x 242       & 0.003125   & 1332140	& 2.877		& 16.211	& 203.035	& 120.117	& 131.252	\\ \hline
			893 x 1279 x 482      & 0.0015625  & 5346482	& 19.807	& 126.331	& 1578.07	& 1441.91	& 1394.25	\\ \hline
	\end{tabular}}
	\caption{CPU time comparison for distance function calculation methods tested on Teddy Bear triangulated surface data. CPU time was measured in seconds.}
	\label{tab:bear_tm_cpu_time}
\end{table}

\clearpage

%%%%%%%%%%%%%%%%%%%%%%%%%%%%%%%%%%%%%%%%%%%%%%%%%%%Searching for the middle surface%%%%%%%%%%%%%%%%%%%%%%%%%%%%%%%%%%%%%%%%%%%%%%%%%%%%%

\section{Numerical methods for computing the middle surface}\label{search_middle_surface}

While analyzing the algorithms for distance function calculation, we realized that methods that track the source of the distance, such as VDT and DP, can be straightforwardly modified for the search of middle surfaces between data sets. In fact, the main inspiration for us was DP method where we expected such modification should work. We propose how to adjust all previously described algorithms to find the middle surface for more data sets of various kinds. Our approach is based on information propagation throughout which we track the source of the information.
\par In the pseudo-codes Alg. \ref{Initialization_pc} and Alg. \ref{Initialization_tm} we showed how to initialize the distance function from one data set. When we have more data sets, we apply one of the algorithms for them separately on the same computational grid. A change which needs to be applied is that in the array $s_{i,j,k}$ for the source of $d_{i,j,k}$ calculation, we need to track also the information to which data set this source belongs to.
\par The change in the VDT and DP algorithms for our new purpose is very easy because they already include source tracking. Again, what we need to change is to track also a label of the data set from which the information propagates. The modification of the FSM and FMM algorithms is not trivial. These methods originally do not contain any information about sources, thus we need to include it in a proper manner.
\par We display the modification of the FSM algorithm in the pseudo-code Alg. \ref{Modified_FSM_pseudo_code}. In every iteration of the algorithm when we cycle through the grid points, we take the distance value from the neighboring points to solve a quadratic equation. We need to keep track, from which neighbors the distance values enter the quadratic equation, thus we save the indexes $\left(r,s,t\right)$, $r,s,t\in\left\{-1,0,1\right\}$, which identify them. When the solution is calculated for the equation we add up the indexes $\left(r,s,t\right)$, see line 20 of Alg. \ref{Modified_FSM_pseudo_code}, and this will show us which source to save for the current grid point from the sources of its $26$ neighbors.
\par \textcolor{black}{For the FMM algorithm, the modification is shown in the pseudo-code Alg. \ref{Modified_FMM_pseudo_code}. In this modification after a nonfixed voxel is tagged as 'visited' (the visiting value is set to $2$), we analyze all its neighbors, see line 26 of the pseudo-code. With the neighbors that also were 'visited', we calculate the current voxel's possible distance from the neighbors' sources, which for quick calculation will be determined by the neighbors' distance value plus the distance between the voxel and its neighbor. The source for which the calculated value is the smallest will be chosen as the source for the current voxel.}
\par By using the modified algorithms, the grid points in the computational domain will be divided into subvolumes "belonging" to the different data sets by source information propagation. To obtain the middle surface between the data sets we just need to find the borders between these subvolumes. To that goal, we use two methods. For any number of data sets, we can cycle through all points of the computational domain and find every point which has a neighbor belonging to a different subvolume. \textcolor{black}{If we apply this for every data set separately, for each of them we obtain a set of points which are at a discrete border of the subvolume belonging to it.} If we have just two data sets, we can treat the obtained information about which data set the grid points belong to, as a function of values $0$ or $1$, and visualize the isosurface of the function with the value $0.5$. \textcolor{black}{We demonstrate the two approaches of visualizing the results in the next subsection with the first numerical experiment for finding the middle surface. In Figure \ref{fig:SpongeSphere_FSM_experiment} in the second picture of the right column we see the representation of the middle surface as a discrete border of subvolumes belonging to a data set, and in the third picture of the right column as an isosurface of a function.}

\begin{algorithm}
	\caption{Modified fast sweeping method including sources}
	\label{Modified_FSM_pseudo_code}
	\begin{algorithmic}[1]
		\renewcommand{\algorithmicrequire}{\textbf{Input:}}
		\Require \texttt{From correct initialization: 3D grid, $d_{i,j,k}$, $f_{i,j,k}$, $s_{i,j,k}$ }
		\renewcommand{\algorithmicrequire}{\textbf{Definition:}}
		%\Require \texttt{$8$ alternating orderings as in Alg.\ref{FSM_pseudo_code}}
		\For{$\left(l=0;l<8;l=l+1\right)$}
		\For{$\left(i=i_{sweep}\left[l,0\right];i \leq i_{sweep}\left[l,1\right];i=i+i_{sweep}\left[l,2\right]\right)$}
		\For{$\left(j=j_{sweep}\left[l,0\right];j \leq j_{sweep}\left[l,1\right];j=j+j_{sweep}\left[l,2\right]\right)$}
		\For{$\left(k=k_{sweep}\left[l,0\right];k \leq k_{sweep}\left[l,1\right];k=k+k_{sweep}\left[l,2\right]\right)$}
		\If{$f_{i,j,k}$ \texttt{is not true} }
		\State \texttt{The indexes $\left(r,s,t\right)$, $r,s,t\in\left\{-1,0,1\right\}$, indicate from}
		\State \texttt{which neighbor the distance value comes from.}
		\State $\left[a_{1},\left(r,s,t\right)_{a_{1}}\right] = \underset{d}{min} \left(\left[d_{i+1,j,k},\left(1,0,0\right)\right], \left[d_{i-1,j,k},\left(-1,0,0\right)\right]\right)$
		\State $\left[a_{2},\left(r,s,t\right)_{a_{2}}\right] = \underset{d}{min} \left(\left[d_{i,j+1,k},\left(0,1,0\right)\right], \left[d_{i,j-1,k},\left(0,-1,0\right)\right]\right)$
		\State $\left[a_{3},\left(r,s,t\right)_{a_{3}}\right] = \underset{d}{min} \left(\left[d_{i,j,k+1},\left(0,0,1\right)\right], \left[d_{i,j,k-1},\left(0,0,-1\right)\right]\right)$ 
		\State \Comment{Use $+\infty$ if $\left(i,j,k\right)$ is out of bounds.}
		\State \texttt{Sort $\left \{ \left[a_{1},\left(r,s,t\right)\right], \left[a_{2},\left(r,s,t\right)\right], \left[a_{3},\left(r,s,t\right)\right] \right \}$ from lowest to} 
		\State \texttt{highest according to values $\left \{ a_{1}, a_{2}, a_{3} \right \}$.}
		\State $\left[d_{new},\left(r,s,t\right)_{d_{new}}\right]=\left[a_{1},\left(r,s,t\right)_{a_{1}}\right]+\left[h,\left(0,0,0\right)\right]$
		\If{$d_{new}>a_{2}$}
		\State $d_{new} = \underset{x}{MaxSolution} \left( \left( x - a_{1}  \right)^{2} + \left( x - a_{2}  \right)^{2} = h^{2} \right)$
		\State $\left(r,s,t\right)_{d_{new}} = \left(0,0,0\right)+\left(r,s,t\right)_{a_{1}}+\left(r,s,t\right)_{a_{2}}$
		\If{$d_{new}>a_{3}$}
		\State $d_{new} = \underset{x}{MaxSolution} \left( \left( x-a_{1}  \right)^{2} + \left( x-a_{2}  \right)^{2} + \left( x-a_{3}  \right)^{2} = h^{2} \right)$
		\State $\left(r,s,t\right)_{d_{new}} = \left(0,0,0\right)+\left(r,s,t\right)_{a_{1}}+\left(r,s,t\right)_{a_{2}}+\left(r,s,t\right)_{a_{3}}$
		\EndIf
		\EndIf
		%\If{$d_{new}<d_{i,j,k}$}
		%\State $d_{i,j,k}=d_{new}$
		%\State $s_{i,j,k}=s_{\left(i,j,k\right)+\left(r,s,t\right)_{d_{new}}}$
		%\EndIf
		\State \textbf{if} \texttt{$d_{new}<d_{i,j,k}$} \textbf{then} \texttt{$\left\lbrace d_{i,j,k}=d_{new}; \; s_{i,j,k}=s_{\left(i,j,k\right)+\left(r,s,t\right)_{d_{new}}} \right\rbrace$}
		\EndIf
		\EndFor
		\EndFor
		\EndFor
		\EndFor
	\end{algorithmic}
\end{algorithm}

\begin{algorithm}
	\caption{Modified fast marching method including sources}
	\label{Modified_FMM_pseudo_code}
	\begin{algorithmic}[1]
		\renewcommand{\algorithmicrequire}{\textbf{Input:}}
		\Require \texttt{From Alg.\ref{Initialization_pc}: 3D grid, $d_{i,j,k}$, $f_{i,j,k}$, $s_{i,j,k}$ }
		%\renewcommand{\algorithmicrequire}{\textbf{Definition:}}
		%\Require \texttt{Set of 6 neighbors $P^{1}$ \eqref{set:6_neighbors} and set of 26 neighbors $P^{2}$ \eqref{set:26_neighbors}.}
		\renewcommand{\algorithmicrequire}{\textbf{Declaration:}}
		\Require \texttt{$v_{i,j,k}$ will hold the visiting values of grid points} \linebreak
		\texttt{'unvisited'=0, 'to be visited'=1, 'visited'=2}
		\Require \texttt{$heap$ container will be a \textit{min-priority-heap}}
		\renewcommand{\algorithmicrequire}{\textbf{Initialization:}}
		\Require \texttt{$\forall f_{i,j,k} = true: \left\lbrace v_{i,j,k} = 1; \; heap.InsertNode(d_{i,j,k}) \right\rbrace $ else: $v_{i,j,k} = 0$}
		\While{\texttt{$heap$ is not empty}}
			\State $\left(i,j,k\right) = heap.GetRoot()$ \Comment{Obtain $(i,j,k)$ with minimum $d$ and delete from $heap$.}
			\For{\texttt{all $\left \{ \left( i + r, j + s, k + t \right) ; \left(r,s,t\right) \in P^{1}\right \}$, not out of bound}}
				\Statex \textbf{.}
				\Statex \textbf{.} \Comment{The pseudo-code is the same as in Alg. \ref{FMM_pseudo_code}.}
				\Statex \textbf{.}
				\setcounter{ALG@line}{21}
			\EndFor
			\State $v_{i,j,k} = 2$
			\If{$f_{i,j,k}$ \texttt{is not true} }
				\State $\left[d_{min},\left(u,v,w\right)\right] = \left[\infty,\left(0,0,0\right)\right]$
				\For{\texttt{all $\left \{ \left( i + r, j + s, k + t \right) ; \left(r,s,t\right) \in P^{2}\right \}$, not out of bound}}
					\If{\texttt{$v_{i+r,j+s,k+t} = 2$}}
						\State \textbf{if} $\left|r\right|+\left|s\right|+\left|t\right|=1$ \textbf{then} $d_{test} = d_{i + r,j + s,k + t} + h$
						\State \textbf{if} $\left|r\right|+\left|s\right|+\left|t\right|=2$ \textbf{then} $d_{test} = d_{i + r,j + s,k + t} + \sqrt{2}*h$
						\State \textbf{if} $\left|r\right|+\left|s\right|+\left|t\right|=3$ \textbf{then} $d_{test} = d_{i + r,j + s,k + t} + \sqrt{3}*h$
						\If{$d_{min}>d_{test}$}
							\State $\left[d_{min},\left(u,v,w\right)\right] = \left[ d_{test},\left( i + r, j + s, k + t \right) \right]$
						\EndIf
					\EndIf				
				\EndFor
				\State $s_{i,j,k}=s_{\left(u,v,w\right)}$
			\EndIf
		\EndWhile
	\end{algorithmic}
\end{algorithm}

\subsubsection{Experiment 1: Sponge \& Sphere}\label{sec:experiment1_sponge_sphere}

\begin{figure}[htp]
	\centering
	\includegraphics[width=\linewidth]{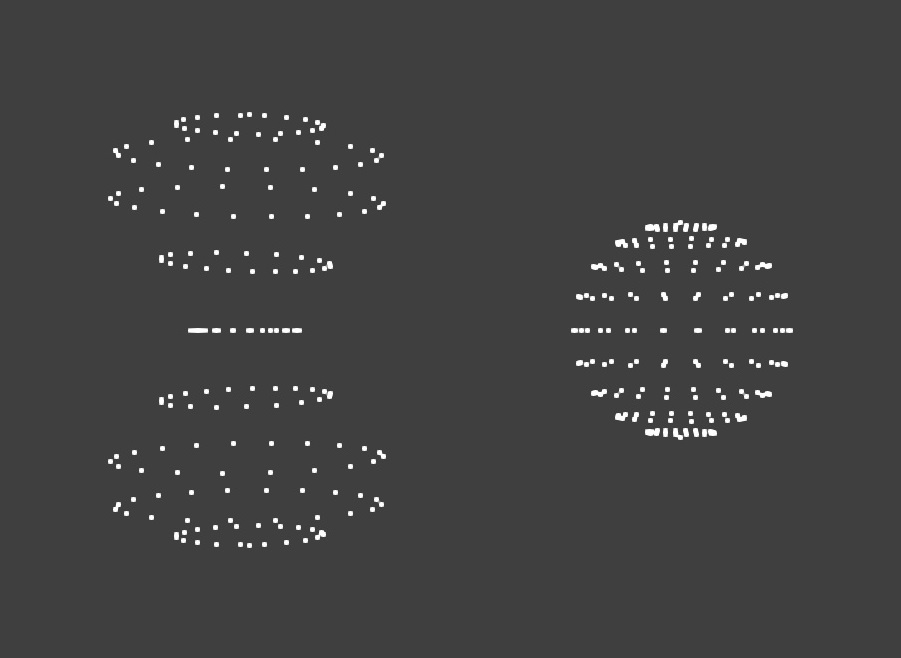}
	\caption{Experiment 1: Generated point cloud data of Sponge and Sphere.}
	\label{fig:SpongeSphere_pc}
\end{figure}

\begin{figure}[htp]
	\centering
	\includegraphics[width=\linewidth]{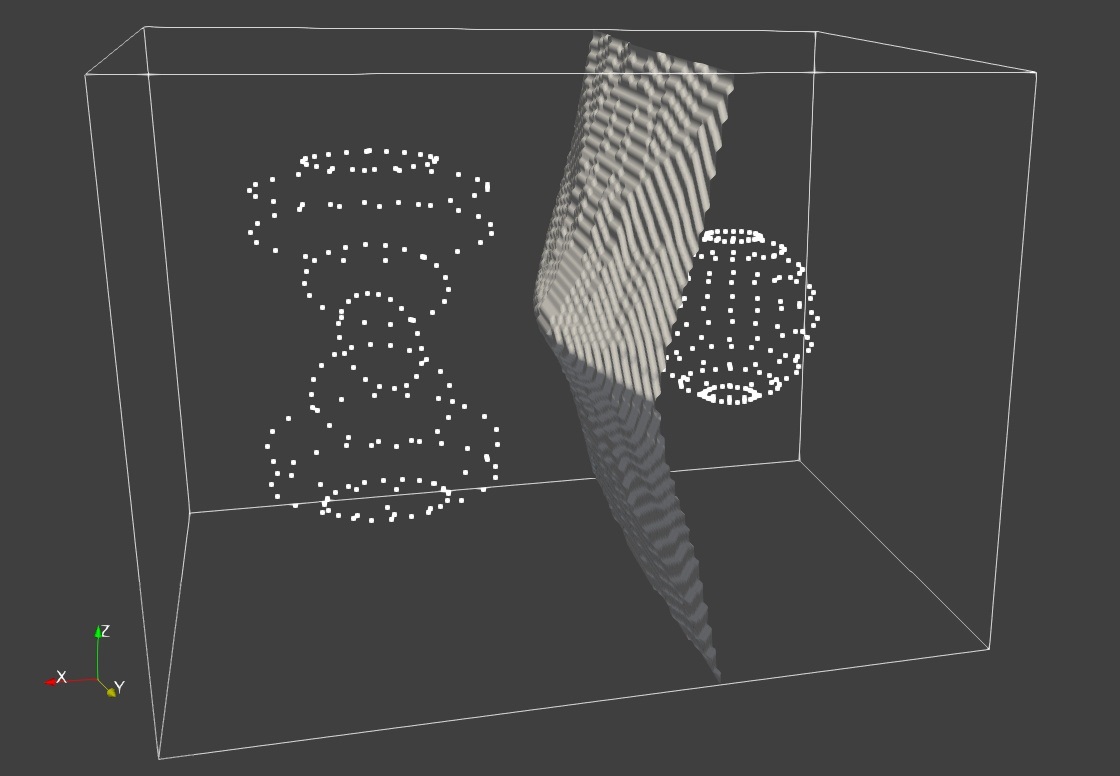}
	\caption{Experiment 1: Middle surface between Sponge and Sphere point cloud data calculated by the VDT algorithm.}
	\label{fig:SpongeSphere_VDT_result}
\end{figure}

Let us have two 3D point clouds, presented in Figure \ref{fig:SpongeSphere_pc}. The first is the "Sponge" point cloud data created by the parametric equations 
\begin{equation} \label{equ:Sponge_parametric}
	\begin{split}
		&x=s_{x} + \left( 0.207 + 2.003 \cdot sin^{2}\left(\varphi\right) - 1.123 \cdot sin^{4}\left(\varphi\right) \right) \cdot cos\left(\varphi\right) \cdot sin\left(\theta\right),\\
		&y=s_{y} + cos\left(\varphi\right) \cdot sin\left(\theta\right),\\
		&z=s_{z} + sin\left(\varphi\right),\\
		&\varphi \in \left \langle 0,2\pi \right ), \theta \in \left \langle 0,\pi \right ).
	\end{split}
\end{equation}
The second point cloud data is a sphere with a radius of $0.5$. The distance between the centers of the two objects is $2.0$. To create the point cloud data we used a step of $\frac{\pi}{10}$ for both angles in the parametric equations. We calculated the distance function on a grid voxel edge size $0.025$. For the middle surface calculated by the VDT algorithm, we obtained the result seen in Figure \ref{fig:SpongeSphere_VDT_result} visualized as an isosurface. The FMM and DP algorithms yield a similar result.

\par
\begin{figure}[htp]
	\centering
	\begin{minipage}{0.49\linewidth}
		\centering
		\includegraphics[width=\linewidth]{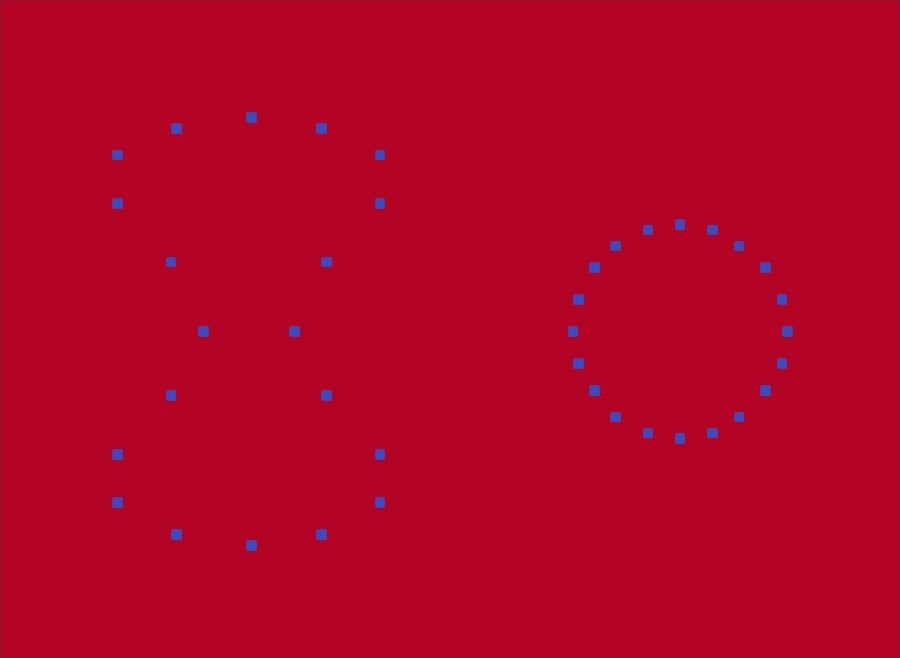}
	\end{minipage}
	\vspace{0.08cm}
	%\hspace{0.04cm}
	\begin{minipage}{0.49\linewidth}
		\centering
		\includegraphics[width=\linewidth]{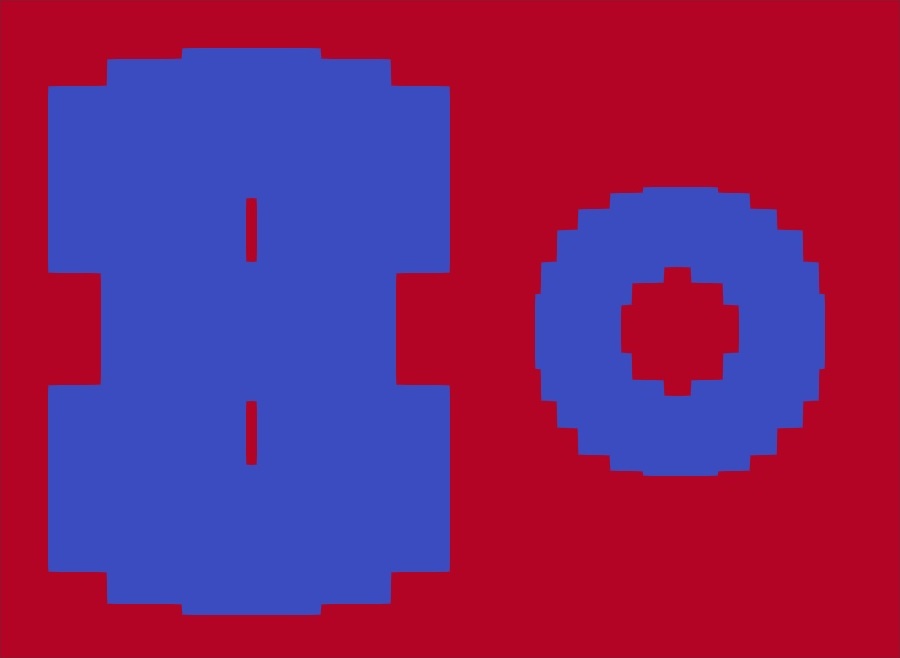}
	\end{minipage}
	\vspace{0.08cm}
	\begin{minipage}{0.49\linewidth}
		\centering
		\includegraphics[width=\linewidth]{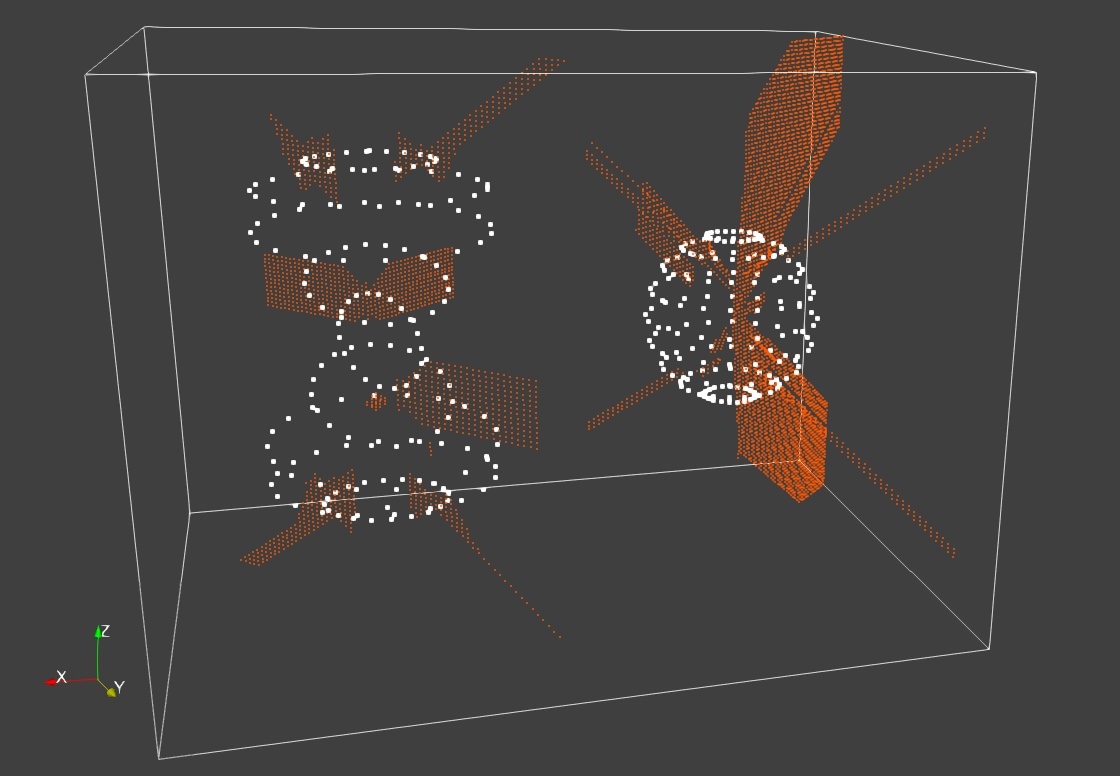}
	\end{minipage}
	\vspace{0.08cm}
	\begin{minipage}{0.49\linewidth}
		\centering
		\includegraphics[width=\linewidth]{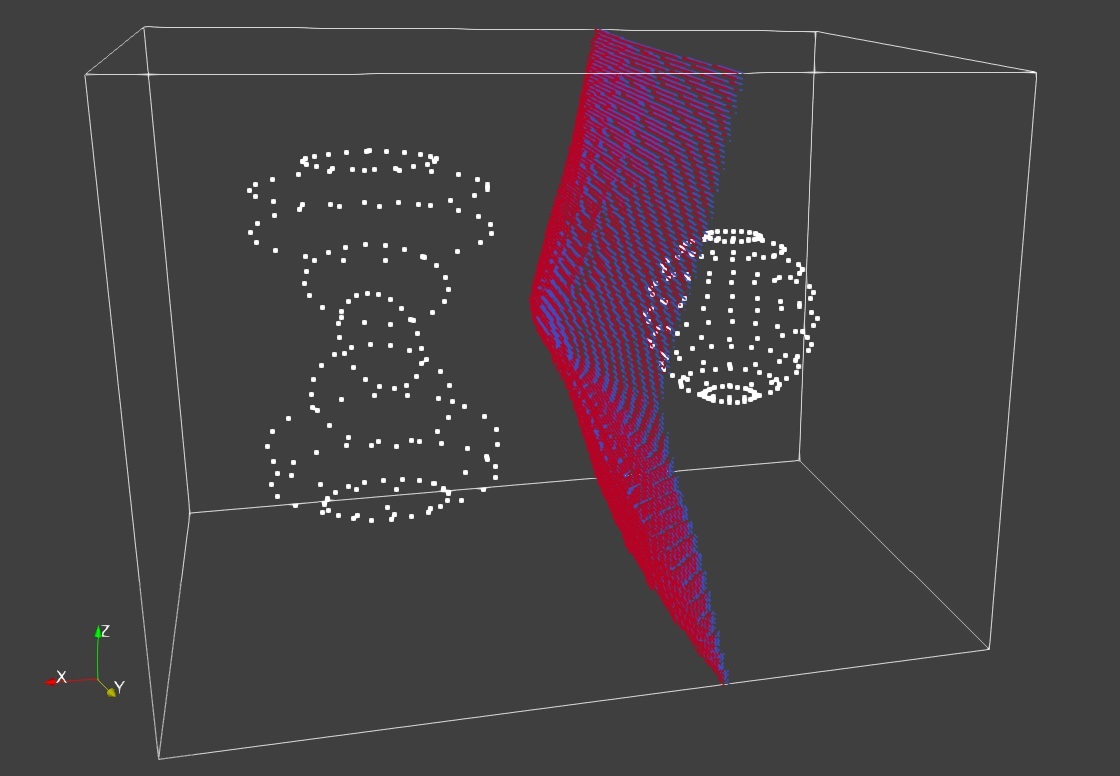}
	\end{minipage}
	\vspace{0.08cm}
	\begin{minipage}{0.49\linewidth}
		\centering
		\includegraphics[width=\linewidth]{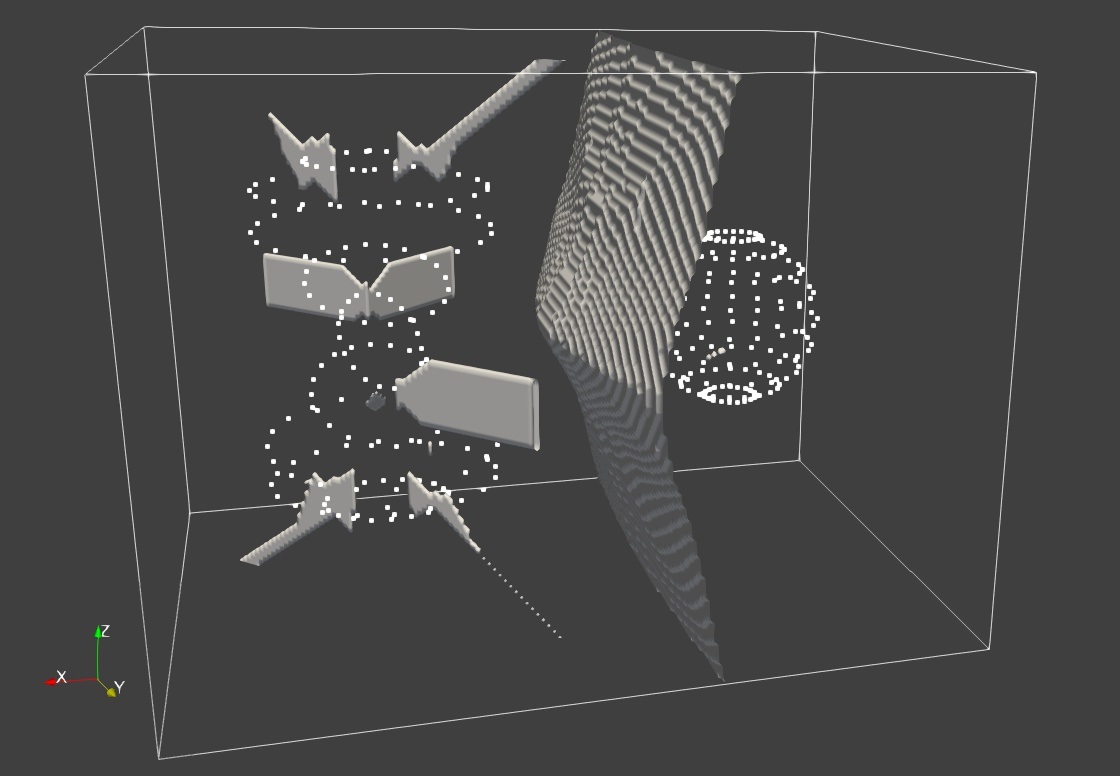}
	\end{minipage}
	%\vspace{0.08cm}
	\begin{minipage}{0.49\linewidth}
		\centering
		\includegraphics[width=\linewidth]{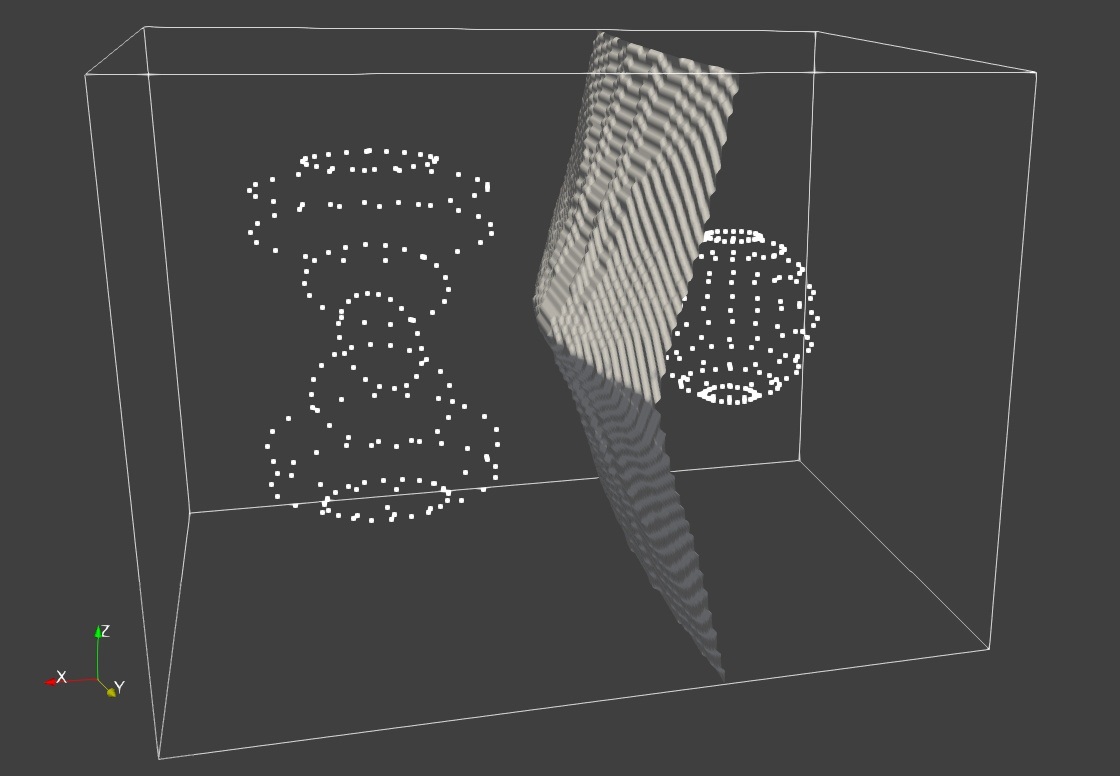}
	\end{minipage}
	%\vspace{0.08cm}
	\caption{Experiment 1: Finding the middle surface between Sponge and Sphere point cloud data by the FSM algorithm. In the first picture of the left column, we see the section of the original initial condition in a constant $y$ plane. In the second picture of the left column, the grid points with no source are visualized. In the third picture of the left column, the incorrect isosurface between subvolumes of the computational grid is visualized. In the first picture of the right column, we can see the corrected initial condition. In the second picture of the right column, the discrete borders of the subvolumes obtained by the corrected calculation are visualized. Red points "belong" to Sponge point cloud data blue points "belong" to Sphere point cloud data. In the third picture of the right column, the correct isosurface between subvolumes of the computational grid is visualized.}
	\label{fig:SpongeSphere_FSM_experiment}
\end{figure}

 With the application of the modified FSM algorithm for this experiment, we discovered that it can cause some issues in specific situations. When we initialize the distance function according to Alg. \ref{Initialization_pc} on a grid with density higher than the point cloud density, we get an initial value that consists of separated subvolumes around the points. The problem is that these gaps in the initialized distance function do not contain any source information, and if we apply FSM, such lack of information can propagate through the computational grid. We can see that in the second picture of the left column in Figure \ref{fig:SpongeSphere_FSM_experiment}. The orange dots indicate the grid points with no source information. This leads to errors when we are trying to detect grid points on the discrete borders of subvolumes belonging to the different point cloud data sets or when we want to visualize the middle surface as an isosurface of a function. The isosurface with errors can be seen in the third picture of the left column in Figure \ref{fig:SpongeSphere_FSM_experiment}.

\par To solve this problem, we need to modify also the initialization of the distance function to point cloud data for the FSM algorithm. The idea is to get a contiguous subvolume for the initialized grid points. For this, we need to increase the volume around the single points in which we initially calculate the distance function. We need to find the minimum size of this volume so that for two neighboring cloud points the volumes will intersect. \textcolor{black}{We found that for this minimum size we can use the maximum of all minimal distances between two cloud points.} With its value, we build a cube around every cloud point which determines the volume in which we will calculate the exact distance values. We can see the result of this modification in the right column of Figure \ref{fig:SpongeSphere_FSM_experiment}. In the first picture, visualizing a section of the new initial condition, we can see that now we have a contiguous subvolume of grid points. In the second picture, we can see that the discrete border of the subvolumes belonging to a data set can be detected correctly, and in the third picture that the isosurface is obtained without any error.
\par In the following experiments, we will show various cases of how we can apply the described algorithms and discuss possible differences in the results of the methods.

\clearpage

\subsubsection{Experiment 2: Subsets of the Cube}\label{sec:experiment2_cube_subsets}

\par
\begin{figure}[htp]
	\centering
	\vspace{1cm}
	\includegraphics[width=\linewidth]{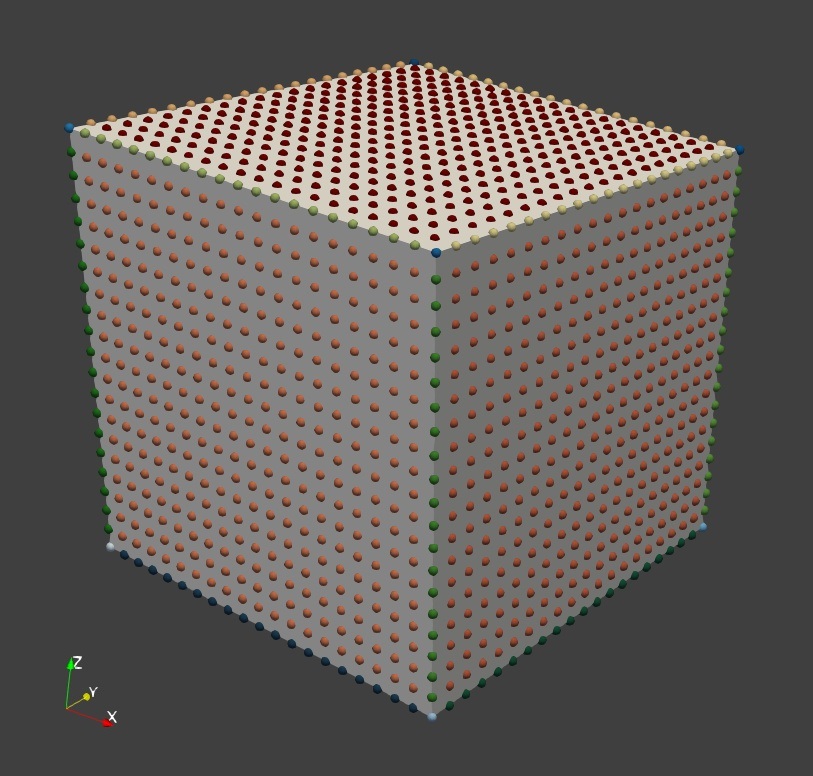}
	\caption{Experiment 2: Visualization of source labels on the Cube data set.}
	\label{fig:cube_every_source}
\end{figure}

 We return to the Cube data set that has coinciding points with the computational grid. We will use a computational grid with voxel edge size $0.05$. The points on every subset of the Cube (vertex, edge, wall) will be treated as a separate data set. In Figure \ref{fig:cube_every_source} we visualize with colors how the points are distributed into sets of sources. We can see that the vertices are treated as one-point data sets, the edges do not contain the vertices and the walls do not contain either the edges or the vertices. Now in this setup, we apply the algorithms for computing the middle surface. 
\par
\begin{figure}[htp]
	\centering
	\includegraphics[width=0.87\linewidth]{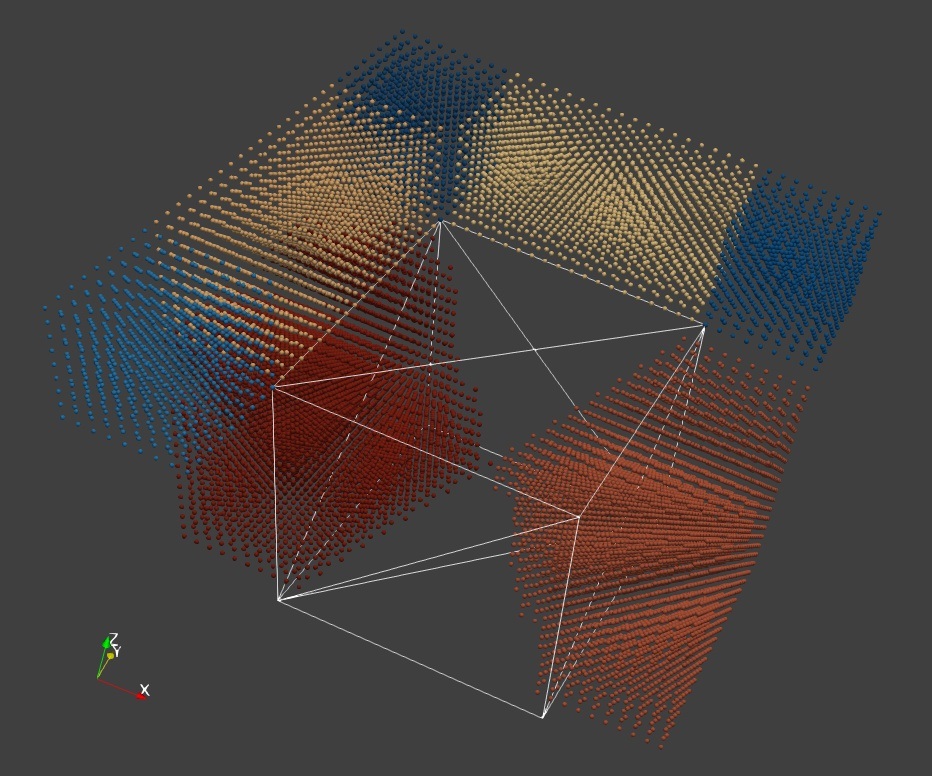}
	\caption{Experiment 2: Visualization of source tracking result on the Cube data set for the VDT method.}
	\label{fig:cube_source_result_VDT}
\end{figure}

\begin{figure}[htp]
	\centering
	\includegraphics[width=0.87\linewidth]{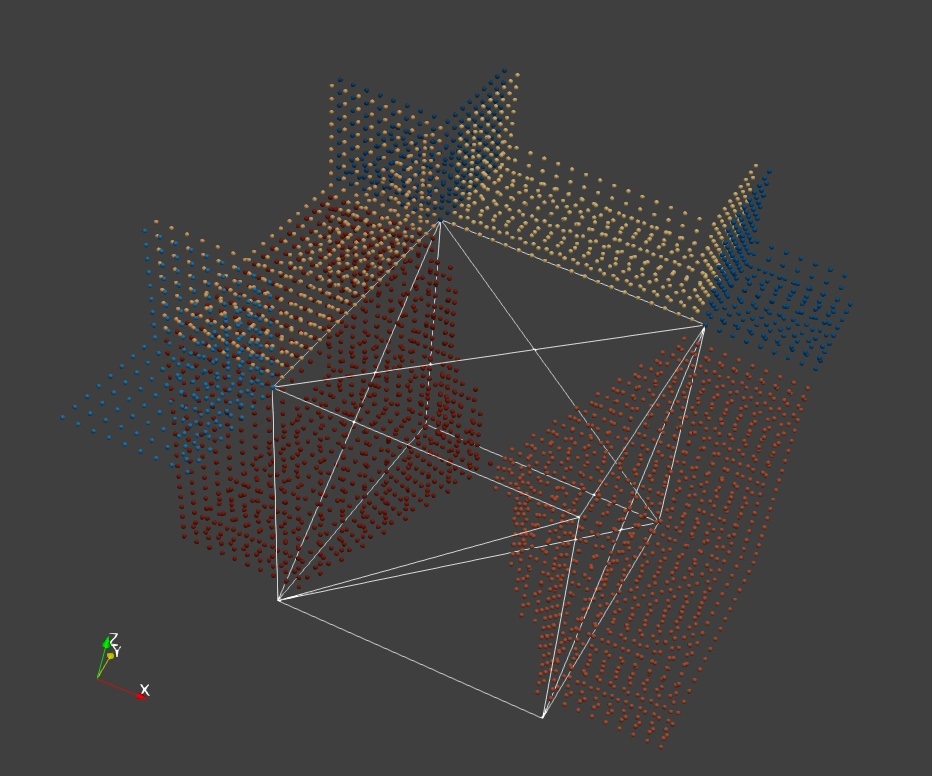}
	\caption{Experiment 2: Visualization of discrete borders of subvolumes belonging to different sources on the Cube data set for the VDT method.}
	\label{fig:cube_source_border_result_VDT}
\end{figure}

 First, we analyze the results for VDT. In Figure \ref{fig:cube_source_result_VDT} we can see how the grid points are assigned to the different subsets. For clearer visualization, we show just some of the separate volumes with the outlines of the Cube by white lines. We can identify by color to which subset of the Cube the points belong to. Let us notice that to the interior of the Cube only the information from the walls propagates. From the vertices and edges, the information only propagates outwards. For this reason, the discrete borders of the subvolumes inside of the Cube are not "uniform".  We can see it more clearly in Figure \ref{fig:cube_source_border_result_VDT} where we visualize only the borders of the subvolumes. For FMM and DP we obtain similar results.
\par
\begin{figure}[htp]
	\centering
	\includegraphics[width=0.87\linewidth]{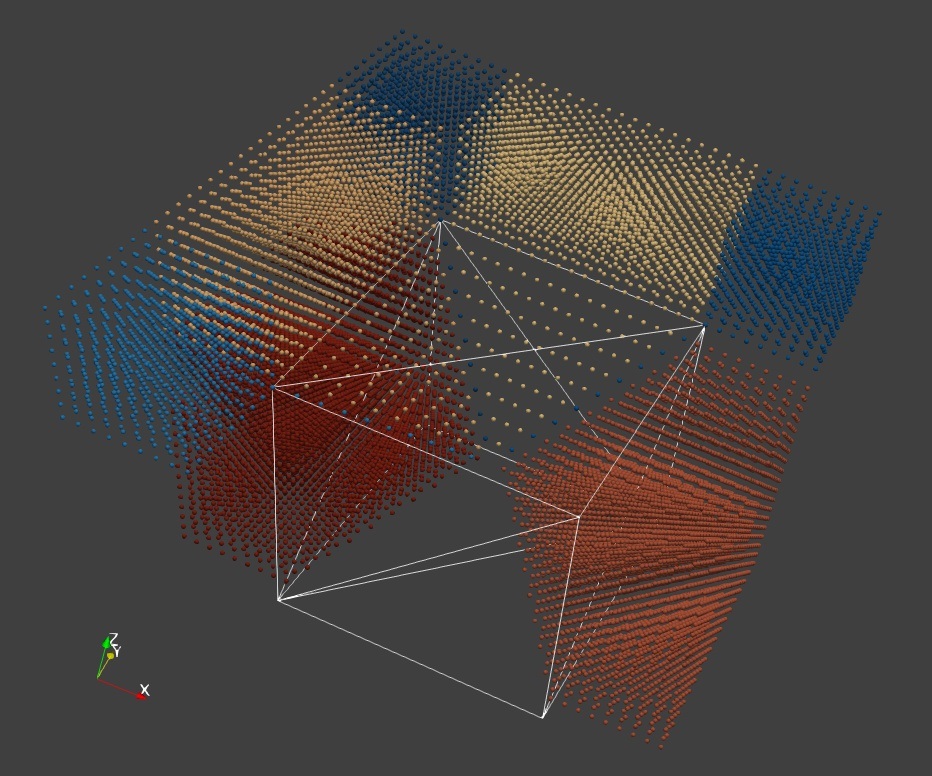}
	\caption{Experiment 2: Visualization of source tracking result on the Cube data set for the FSM algorithm.}
	\label{fig:cube_source_result_FSM}
\end{figure}

\begin{figure}[htp]
	\centering
	\includegraphics[width=0.87\linewidth]{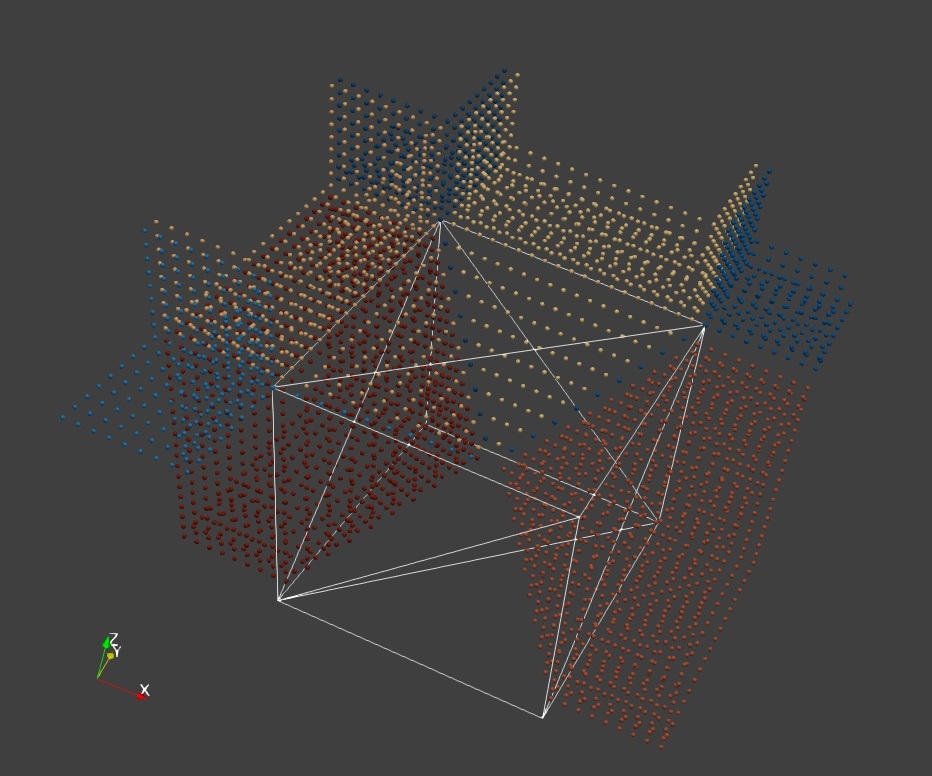}
	\caption{Experiment 2: Visualization of discrete borders of subvolumes belonging to different sources on the Cube data set for the FSM algorithm.}
	\label{fig:cube_source_border_result_FSM}
\end{figure}

 Let us compare the previous result to the results of the FSM algorithm visualized in Figure \ref{fig:cube_source_result_FSM}. We can identify by the colors that the information propagates inward from all subsets of the Cube. Inside of the Cube, the grid points belonging to vertices are along a line, for edges, the grid points are confined to a triangle, and for the walls, they are inside a pyramid. In these results, the borders of the separated volumes are much clearer and sharper which we can identify easier in Figure \ref{fig:cube_source_border_result_FSM}.

\clearpage

\subsubsection{Experiment 3: Cube \& Sphere}\label{sec:experiment3_cube_sphere}

\begin{figure}[htp]
	\centering
	\includegraphics[width=0.49\linewidth]{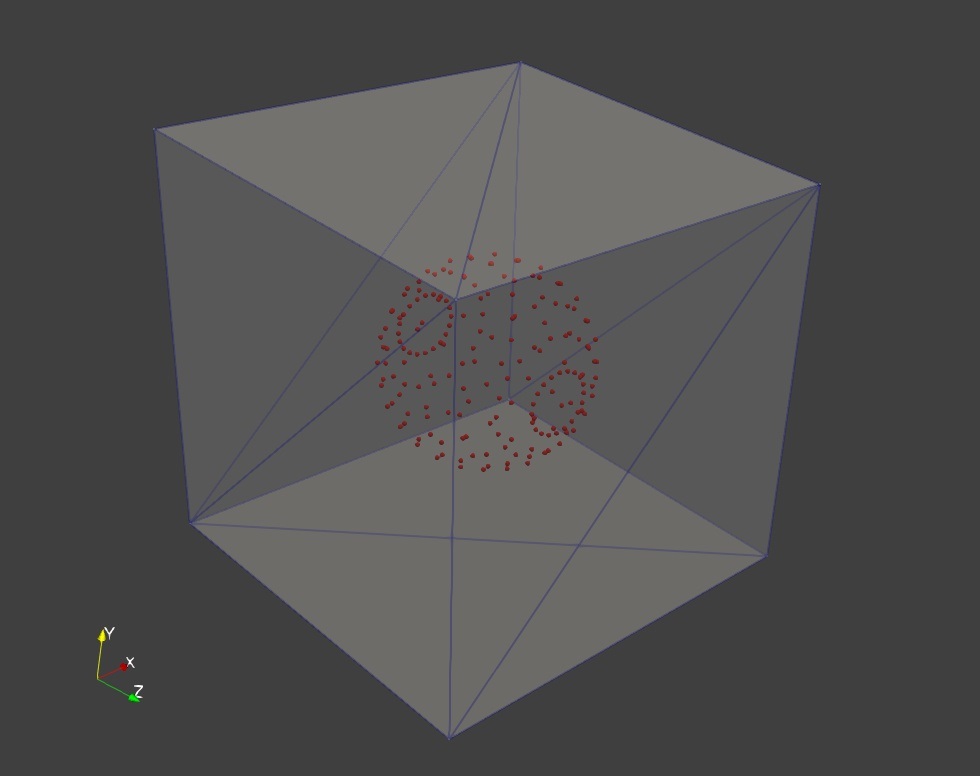}
	\includegraphics[width=0.49\linewidth]{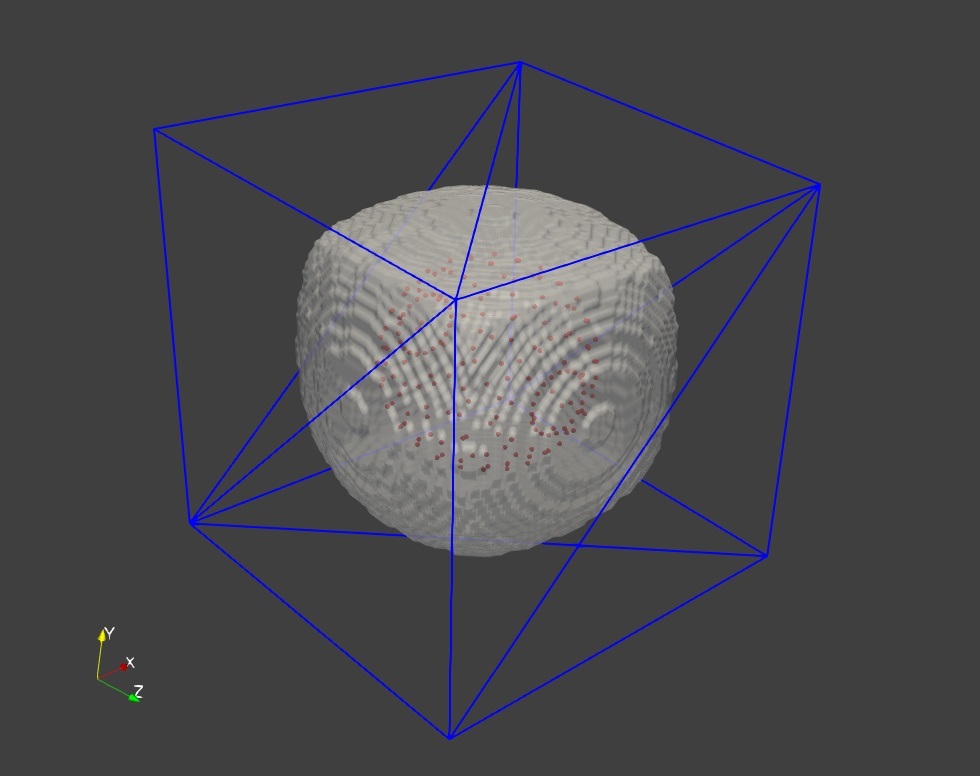}
	
	\vspace{0.08cm}
	
	\includegraphics[width=0.49\linewidth]{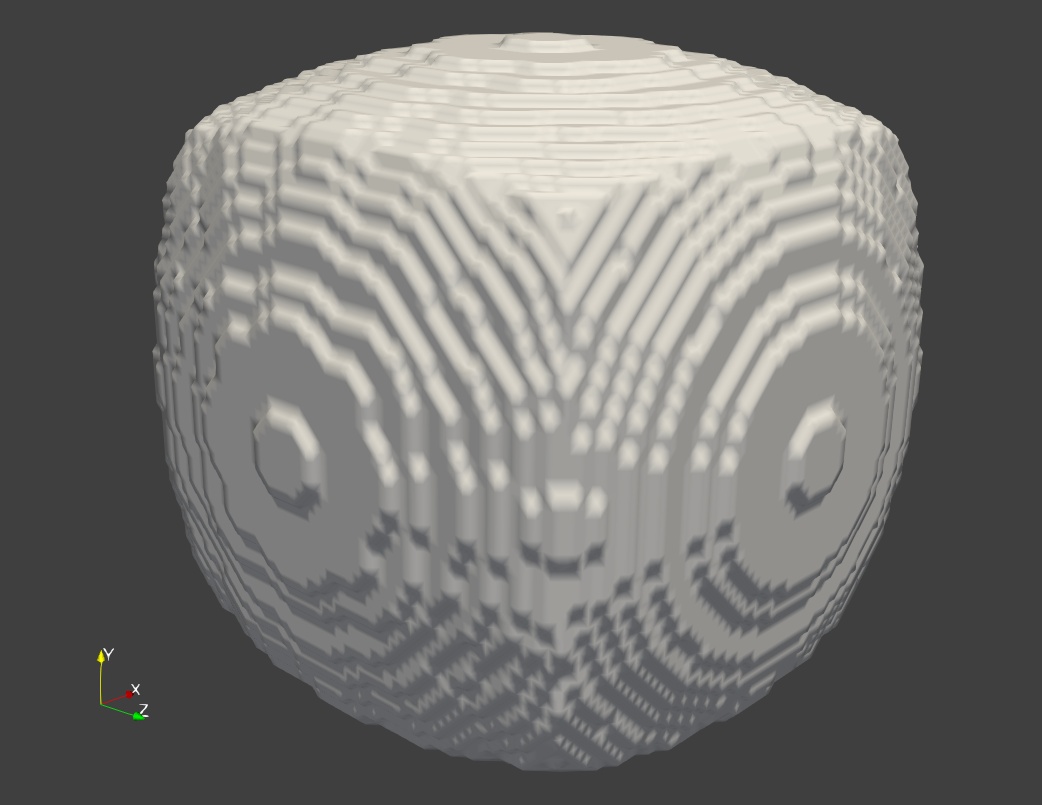}
	\includegraphics[width=0.49\linewidth]{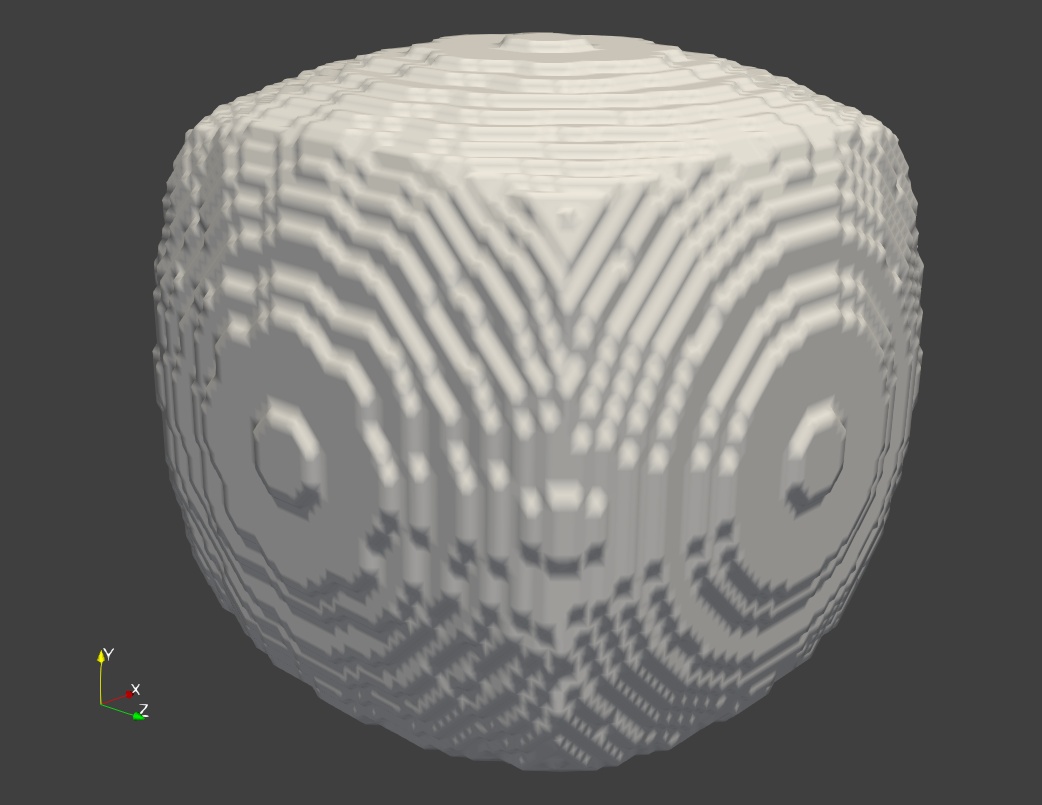}	
	
	\vspace{0.08cm}
	
	\includegraphics[width=0.49\linewidth]{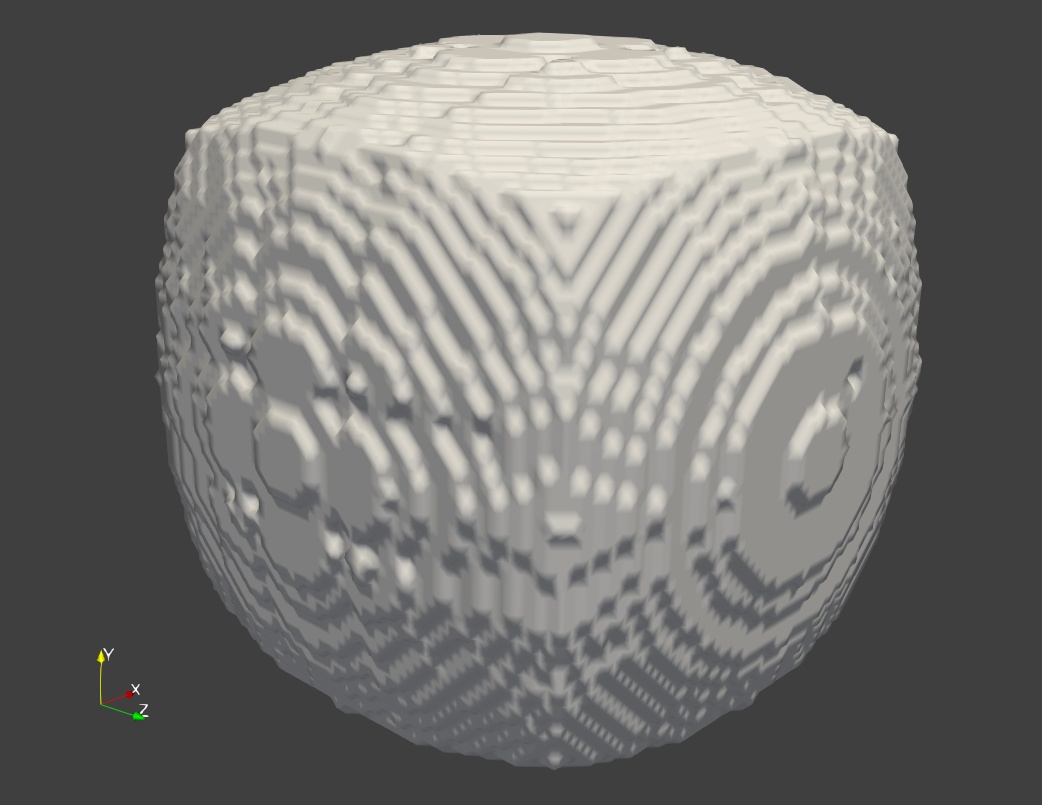}
	\includegraphics[width=0.49\linewidth]{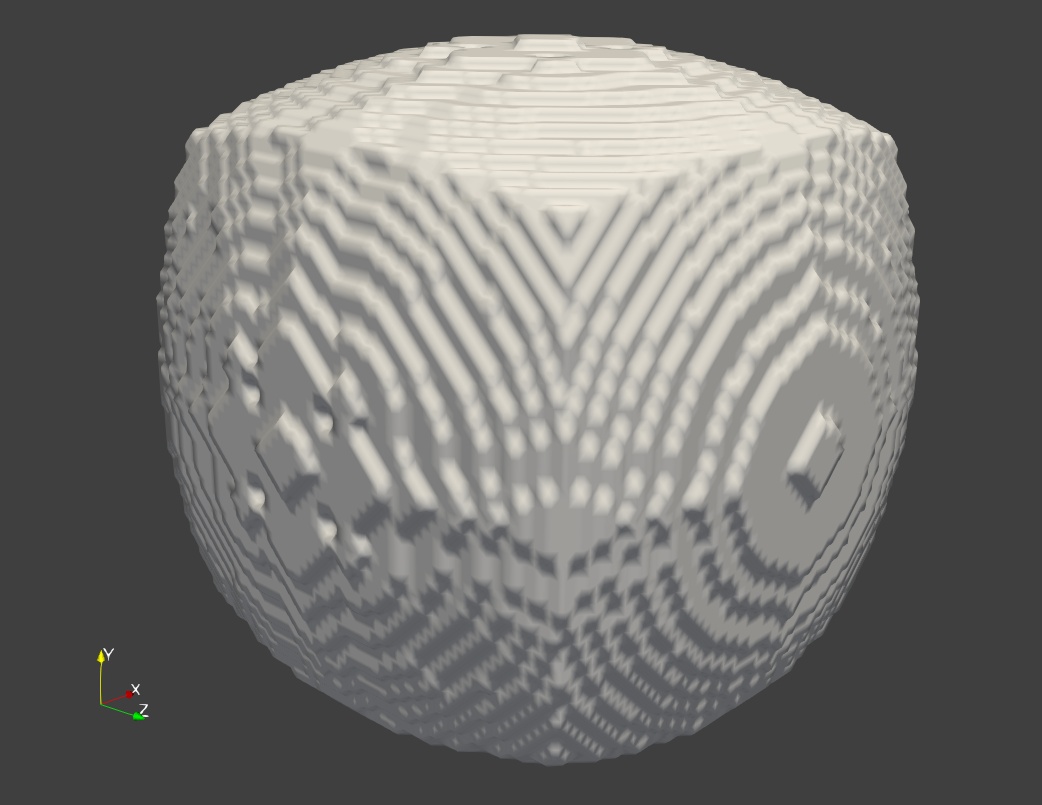}
	
	\caption{Experiment 3: Finding the middle surface between Sphere point cloud data inside a Cube triangulated surface. In the first picture, we see the two objects. In the next picture, we visualize the middle surface together with the objects. In the following pictures, we show the resulting isosurfaces for every algorithm in more detail. They are visualized from left-up to down-right in the following order: VDT, DP, FSM, FMM. The visualized results were computed on a grid with $181^{3}$ elements and a voxel edge size of $0.01$.}
	\label{fig:CubeSphere_Sphere10_results}
\end{figure}

In the next experiment, we consider a cube with the same parameters but now we will work with it as a triangulated surface. As we demonstrated in Section \ref{sec:distance_to_tm} we can use the algorithms for distance function calculations on triangulated surfaces as well if we use the Alg. \ref{Initialization_tm} for the initialization. This type of initialization produces contiguous subvolumes of grid points thus it does not need any changes to be applicable for the Modified FSM algorithm as it was in the case of point cloud data. Inside of the Cube, we have the Sphere with radius $0.25$ and center point the same as the center of the Cube. We can see their relative location in the first picture of Figure \ref{fig:CubeSphere_Sphere10_results}. In the second picture, we see the computed middle surface with the objects. In the next pictures of this figure, we can see the results for the VDT, DP, FSM, and FMM algorithms, in this order from left-up to right-down. In the detailed view of the results, we can see the fine differences between them.
\par For a quantitative comparison of the methods, we calculate the volume and area of the isosurfaces computed on computational grids with different voxel edge sizes, equal to $0.2$, $0.1$, $0.05$, $0.025$, $0.0125$, $0.00625$, $0.003125$. We list these results in Table \ref{tab:CubeSphere}. By comparing the values in this table and looking at the pictures of the middle surface we can see that the results from the pairs VDT, DP, and FSM, FMM are very similar.

\begin{table}[htp]	
	\centering
	\catcode`\-=12
	\scalebox{0.75}{
	\begin{tabular}{|c|c|c|c|c|c|c|c|c|c|}
		\hline
		Number of   & Voxel     & \multicolumn{2}{c|}{VDT} & \multicolumn{2}{c|}{DP} & \multicolumn{2}{c|}{FSM} & \multicolumn{2}{c|}{FMM} \\ \cline{3-10} 
		grid points & edge size & Volume    & Area    	& Volume    & Area   	& Volume    & Area   	& Volume    & Area    	\\ \hline
		$10^{3}$    & 0.2       & 0.418667  & 2.92008 	& 0.418667	& 2.92008	& 0.418667  & 2.92008   & 0.418667	& 2.92008 	\\ \hline
		$19^{3}$    & 0.1       & 0.3005    & 2.39785 	& 0.2855	& 2.34128	& 0.244167  & 1.98998   & 0.2645	& 2.21841 	\\ \hline
		$37^{3}$    & 0.05      & 0.291396  & 2.40964   & 0.291396	& 2.40964	& 0.271396  & 2.29456   & 0.273396	& 2.31799   \\ \hline
		$73^{3}$    & 0.025     & 0.286294  & 2.39177   & 0.286326	& 2.39452	& 0.277992  & 2.34297   & 0.274508	& 2.34427   \\ \hline
		$145^{3}$   & 0.0125    & 0.287682  & 2.40421   & 0.287686	& 2.40489	& 0.282912  & 2.39768   & 0.277739	& 2.37683   \\ \hline
		$289^{3}$   & 0.00625   & 0.287828  & 2.42233   & 0.287828	& 2.42245	& 0.284826  & 2.40956   & 0.281192	& 2.40664   \\ \hline
		$577^{3}$   & 0.003125  & 0.287996  & 2.42093   & 0.287995	& 2.42094	& 0.286224  & 2.42284   & 0.283977	& 2.41481   \\ \hline
	\end{tabular}}
	\caption{Experiment 3: Comparing volume and area for middle surface between the Cube and the Sphere data sets.}
	\label{tab:CubeSphere}
\end{table} 

\clearpage

\subsubsection{Experiment 4: Five Ellipsoids} \label{sec:experiment4_five_ellipsoids}

\begin{figure}[htp]
	\centering
	\includegraphics[width=\linewidth]{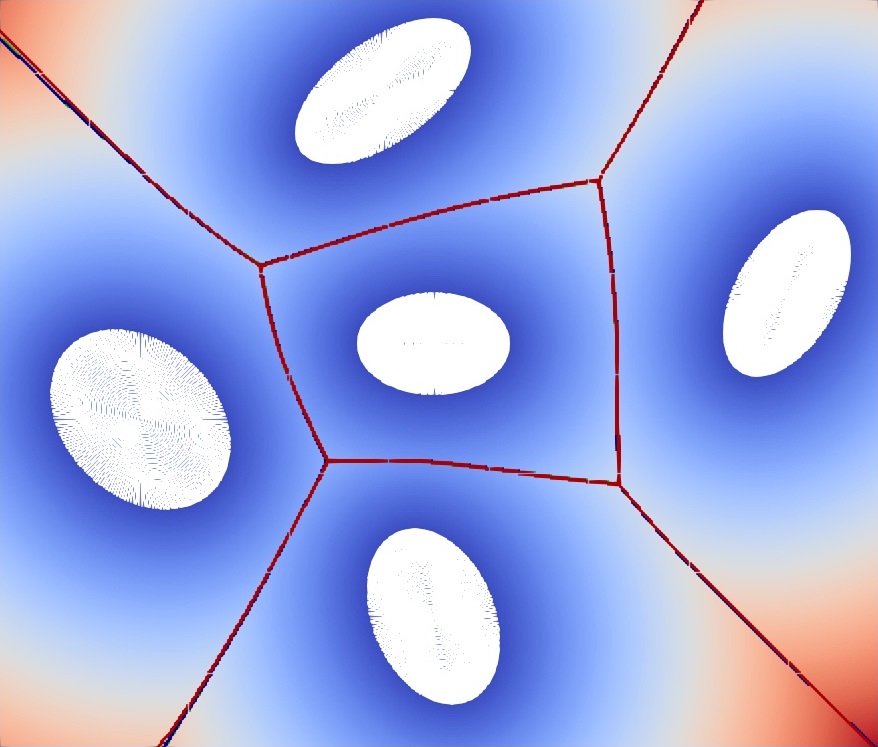}
	\vspace{0.04cm}
	\caption{Experiment 4: \textcolor{black}{Finding the border between five Ellipsoid point cloud data. In the picture, the border points between divided volumes are visualized in the plane $z=0$ together with the distance function and the data sets. We show the results for FSM algorithm with red lines, for VDT with dark blue lines.}}
	\label{fig:FiveEllipsoids_0p020000}
\end{figure}

The following experiment is done with five different Ellipsoids point cloud data sets, for which the center points all lie on the plane $z=0$. In Figure \ref{fig:FiveEllipsoids_0p020000} we visualize the results of the algorithms in the plane $z=0$ as the discrete border of subvolumes together with the distance function and the original data. Here we show the results of the FSM algorithm with red lines, of VDT with dark blue lines. Because of the overlapping of the results for DP and FMM are almost not visible. We can just see the result for DP with a green line in the upper left corner. In this experiment, we can see that with our algorithms the obtained results are a good approximation of the Voronoi diagram.

\clearpage

\subsubsection{Experiment 5: Two parallel surfaces}\label{sec:experiment5_parallel_surfaces}

For the last experiment, we want to show how accurately the algorithms can find the middle surface between two parallel data sets. For this purpose, we will use wave-like surfaces generated as point cloud data by functions 
\begin{equation}\label{eq:cosinus_curves}
	\begin{split}
		&f(x,y)=0.2*cos\left(x*y\right) + 0.5,\\
		&f(x,y)=0.2*cos\left(x*y\right) - 0.5,\\
		&\left(x,y\right)\in<-5.0,5.0>\times <-5.0,5.0>
	\end{split}
\end{equation}
with a step of $0.05$ for both $x$ and $y$ variables. 
\par
\begin{figure}[htp]
	\centering
	\includegraphics[width=0.9\linewidth]{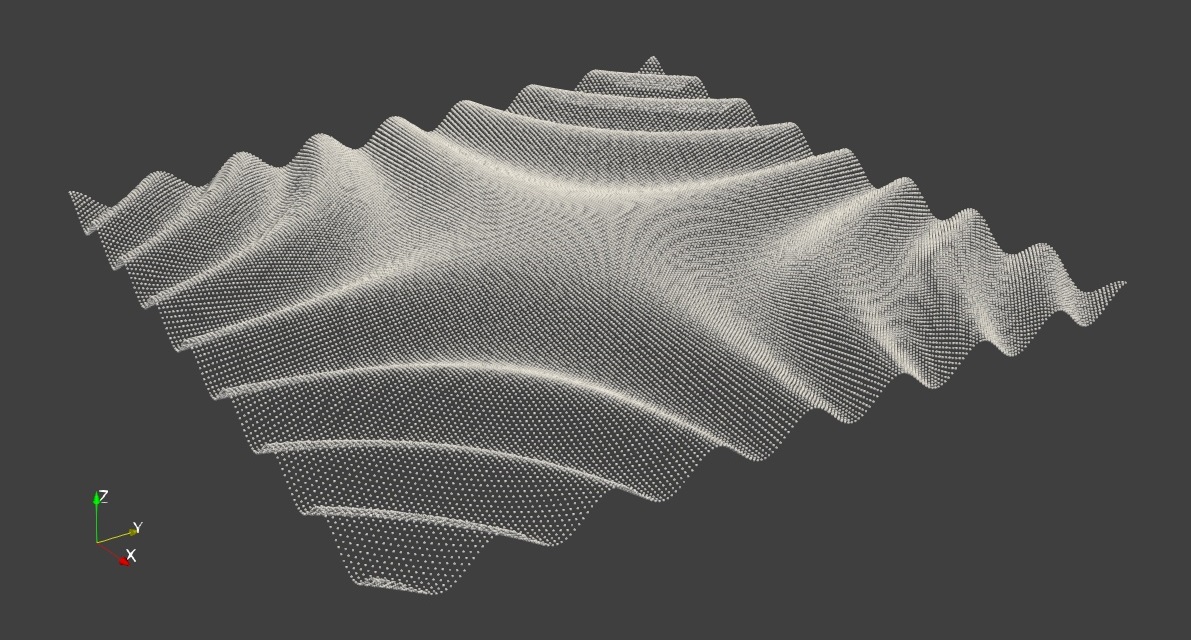}	
	\vspace{0.08cm}
	
	\includegraphics[width=0.9\linewidth]{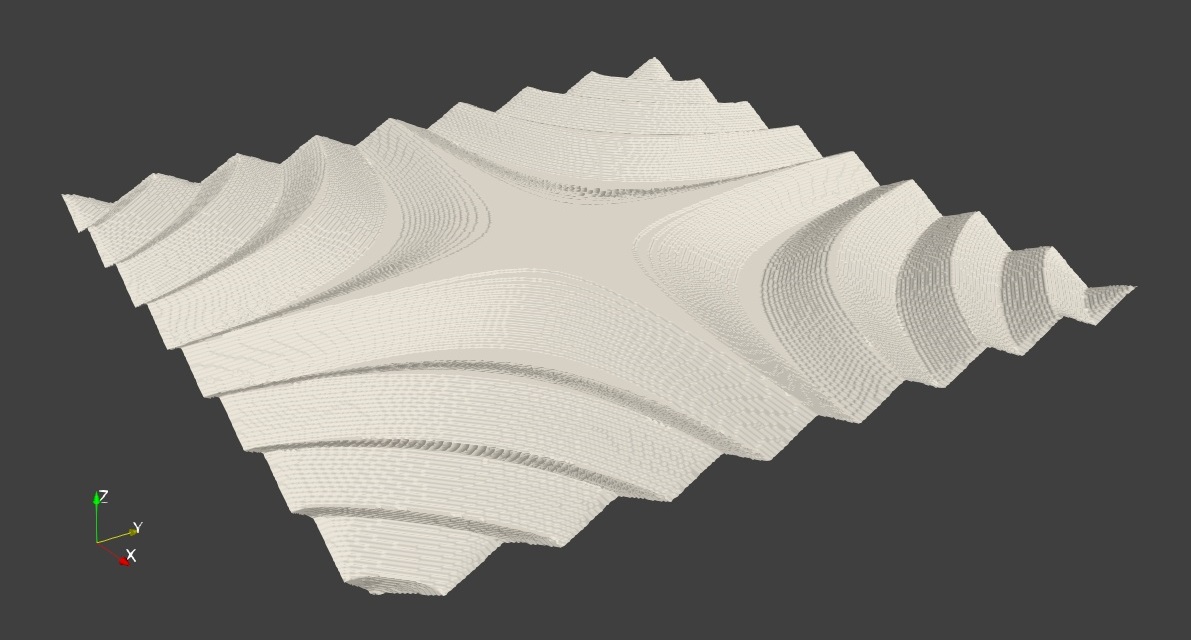}
	\vspace{0.08cm}
	
	\caption{Experiment 5: Finding the middle surface between two parallel wave-like point cloud data sets generated by equations \eqref{eq:cosinus_curves}. In the first picture, we see the visualization of one point cloud. The other one is identical just shifted along the $z$ axis. In the second picture, we visualize the middle surface which divides the computational domain between the two point cloud data sets.}
	\label{fig:CosinusCurves_0p025000}
\end{figure}

 We can see the visualization of the point cloud data generated by the first equation of \eqref{eq:cosinus_curves} in the first picture of Figure \ref{fig:CosinusCurves_0p025000}. In the second picture, we can see the result of the calculations by the FSM algorithm on a computational grid with voxel edge size $0.025$ represented as an isosurface. This isosurface lies between the two parallel point cloud data sets. Visually the results for the four methods do not show noticeable differences, thus we show only the results of FSM.

%%%%%%%%%%%%%%%%%%%%%%%%%%%%%%%%%%%%%%%%%%%%%%%%%%%%%%%%%%%%%%%%%%%%%%%%AcknowledgementAndREFERENCES%%%%%%%%%%%%%%%%%%%%%%%%%%%%%%%%%%%%%%%%%%%%%%%%%%%%%%%%%%%%%%%%%%%%%%%%

\begin{acknowledgements}
We would like to thank Prof. Zuzana Krivá for pointing out the possibility to use 6 voxel neighbors instead of 26 in the Dijkstra-Pythagoras method.
\end{acknowledgements}


\bibliographystyle{spmpsci}      % mathematics and physical sciences
\bibliography{bibfile}   % name your BibTeX data base

\begin{thebibliography}{10}
\providecommand{\url}[1]{{#1}}
\providecommand{\urlprefix}{URL }
\expandafter\ifx\csname urlstyle\endcsname\relax
  \providecommand{\doi}[1]{DOI~\discretionary{}{}{}#1}\else
  \providecommand{\doi}{DOI~\discretionary{}{}{}\begingroup
  \urlstyle{rm}\Url}\fi

\bibitem{Chen_2009_ABF}
Chen, X., Golovinskiy, A., Funkhouser, T.: A benchmark for {3D} mesh
  segmentation.
\newblock ACM Transactions on Graphics \textbf{28}(3) (2009).
\newblock \doi{10.1145/1531326.1531379}

\bibitem{Intro_algo}
Cormen, T.H., Leiserson, C.E., Rivest, R.L., Stein, C.: Introduction to
  algorithms, Third Edition.
\newblock MIT Press and McGraw-Hill (2009)

\bibitem{Danielsson}
Danielsson, P.E.: Euclidean distance mapping.
\newblock Computer Graphics and Image Processing \textbf{14}(3), 227--248
  (1980).
\newblock \doi{10.1016/0146-664X(80)90054-4}

\bibitem{Eberly_PointToTriangle}
Eberly, D.: Distance between point and triangle in {3D}.
\newblock Geometric Tools  (1999).
\newblock
  \urlprefix\url{https://www.geometrictools.com/Documentation/DistancePoint3Triangle3.pdf}

\bibitem{Jones_Baerentzen_Sramek}
Jones, M.W., Baerentzen, J.A., Sramek, M.: {3D} distance fields: a survey of
  techniques and applications.
\newblock IEEE Transactions on Visualization and Computer Graphics
  \textbf{12}(4), 581--599 (2006).
\newblock \doi{10.1109/TVCG.2006.56}

\bibitem{KIMMEL1995382}
Kimmel, R., Shaked, D., Kiryati, N., Bruckstein, A.M.: Skeletonization via
  distance maps and level sets.
\newblock Computer Vision and Image Understanding \textbf{62}(3), 382--391
  (1995).
\newblock \doi{10.1006/cviu.1995.1062}

\bibitem{Persson_phdthesis}
Persson, P.O.: Mesh generation for implicit geometries.
\newblock Ph.D. thesis, Department of Mathematics, Massachusetts Institute Of
  Technology (2005)

\bibitem{Rouy_Tourin}
Rouy, E., Tourin, A.: A viscosity solutions approach to shape-from-shading.
\newblock SIAM Journal on Numerical Analysis \textbf{29}(3), 867--884 (1992).
\newblock \doi{10.1137/0729053}

\bibitem{Rumpf_ContSkelet}
Rumpf, M., Telea, A.: A continuous skeletonization method based on level sets.
\newblock Proceedings of the symposium on Data Visualisation 2002 p. 151–ff
  (2002)

\bibitem{Sethian_FMM}
Sethian, J.A.: A fast marching level set method for monotonically advancing
  fronts.
\newblock Proceedings of the National Academy of Sciences \textbf{93}(4),
  1591--1595 (1996).
\newblock \doi{10.1073/pnas.93.4.1591}

\bibitem{Siddiqi_Hamilton}
Siddiqi, K., Bouix, S., Tannenbaum, A., Zucker, S.: The hamilton-jacobi
  skeleton.
\newblock International Conference on Computer Vision (ICCV)

\bibitem{Smisek_dissertation}
Smíšek, M.: Analysis of {3D} and {4D} images of organisms in embryogenesis.
\newblock Ph.D. thesis, Faculty of Civil Engineering, Slovak University of
  Technology Bratislava (2015)

\bibitem{Models_Archive}
Turk, G., Mullins, B.: Large geometric models archive, {Georgia Institute of
  Technology}  (1999).
\newblock \urlprefix\url{https://www.cc.gatech.edu/projects/large\_models/}

\bibitem{Zhao}
Zhao, H.: A fast sweeping method for {Eikonal} equations.
\newblock Mathematics of Computation \textbf{74}, 603--627 (2005).
\newblock \doi{10.1090/S0025-5718-04-01678-3}

\end{thebibliography}

%\begin{thebibliography}{}
% etc
%\end{thebibliography}

\end{document}